\documentclass{article}

\usepackage{multirow}
\usepackage{graphicx, amssymb, amsfonts, subfigure, setspace, geometry, amsmath}
\usepackage{authblk}
\usepackage{mathtools}

\usepackage[open-square,define-standard-theorems]{QED}
\newtheorem{theorem}{Theorem}[section]

\newtheorem{lemma}{Lemma}[section] 
\newtheorem{definition}{Definition}[section] 
\newtheorem{example}{Example}[section] 
\newtheorem{rem}{Remark}[section]
 
\newcommand{\D}{\ensuremath{\,\mathrm{d}}}
\newcommand{\xb}{\mathbf{x}}
\newcommand{\yb}{\mathbf{y}}

\newcommand{\db}{\mathbf{d}}

\newcommand{\vecp}{\vec{p}_j}

\DeclarePairedDelimiter{\ceil}{\lceil}{\rceil}

\begin{document}
\title{Inverse obstacle scattering in two dimensions with multiple frequency data and multiple angles of incidence}

\author[1]{Carlos Borges}
\author[1,2]{Leslie Greengard}
\affil[1]{Courant Institute of Mathematical Sciences, New York University
              }
\affil[2]{Simons Center for Data Analysis, Simons Foundation, New York, NY}

\date{Received: date / Accepted: date}

\maketitle

\begin{abstract}
We consider the problem of reconstructing the shape of an impenetrable 
sound-soft obstacle from scattering measurements. The input data
is assumed to be the far-field pattern generated when a plane wave 
impinges on an unknown obstacle from one or more directions and at 
one or more frequencies. It is well known that this
inverse scattering problem is both ill posed and nonlinear.  
It is common practice to 
overcome the ill posedness through the use of a penalty method or 
Tikhonov regularization. Here, we present a more physical regularization, 
based simply on restricting the unknown boundary to be band-limited in 
a suitable sense. To overcome the nonlinearity
of the problem, we use a variant of Newton's method.
When multiple frequency data is available, we supplement Newton's 
method with the recursive linearization approach due to Chen.

During the course of solving the inverse problem, we need to compute
the solution to a large number of forward scattering problems. 
For this, we use high-order accurate 
integral equation discretizations, coupled with fast direct solvers when 
the problem is sufficiently large.
\end{abstract}

\section{Introduction}
Inverse problems arise in many parts of science and engineering, 
including medical imaging, remote sensing, ocean acoustics, 
nondestructive testing, geophysics and radar
\cite{Colton,isakov,kirsch_book,kuchment2014radon}. 
In this paper, we concentrate on the problem of recovering the shape of
an unknown obstacle embedded in a homogeneous medium 
from far-field acoustic scattering measurements in 
two space dimensions. We assume 
that the object $\Omega$ is impenetrable and ``sound-soft" with boundary
$\Gamma$. 
We restrict our attention here to the time harmonic setting, in which
case the governing equation is the Helmholtz equation
\[  \Delta u(\xb) + k^2 u(\xb) = 0 \, , \qquad 
\xb \in \mathbb{R}^2 \setminus\overline{\Omega}.
\]

In the sound-soft case,  $u$ must
satisfy the Dirichlet condition
\[ u(\xb) = 0, \qquad {\rm for}\ \xb \in \Gamma.
\]
Following standard practice, we write
\[ u(\xb) = u^{\text{inc}}(\xb) + u^{\text{scat}}(\xb),
\]
where $u^{\text{inc}}$ denotes a known incoming field and 
$u^{\text{scat}}$ denotes the scattered field.
We will assume that the incoming field is a plane wave
$u^{\text{inc}}(\xb)=\exp(i k \xb \cdot \db)$, where 
$\db$ is a unit vector that defines the direction of propagation.

It is well-known that
radiating fields have a simpe asymptotic structure as 
$|\xb| \rightarrow \infty$.
More precisely \cite{Colton},

\begin{lemma} 
Let $(r,\theta)$ denote the polar coordinates of the point $\xb$.
Every radiating solution $u^{\text{\emph{scat}}}$ to the 
Helmholtz equation has the asymptotic behavior of an outgoing 
cylindrical wave:
\begin{eqnarray*}
u^{\text{\emph{scat}}}(\xb)&=&\frac{e^{ik r}}{r^{1/2}}
\left\{u^{\infty}(\theta)+\mathcal{O}\left(\frac{1}{r^{3/2}}\right)\right\},
\quad r \rightarrow \infty,
\end{eqnarray*}
uniformly in all directions $\theta$.
The function $u^{\infty}$ is defined on the unit 
circle $U$ and referred to as the far-field pattern of 
$u^{\text{\emph{scat}}}$. 
\end{lemma}

\begin{definition}
Suppose that $u^{\text{\emph{inc}}}$ and the domain shape $\Gamma$
are known.  The {\em forward} scattering problem consists of 
solving the boundary value problem:
\begin{eqnarray}
\Delta u^{\text{\emph{scat}}}(\xb) +k^2 u^\text{\emph{scat}}(\xb) &=&0 
\quad  \mbox{in} \quad \mathbb{R}^2 \setminus \overline{\Omega}  \label{fscat} \\
u^\text{\emph{scat}}(\xb)&=&  -u^\text{\emph{inc}}(\xb) 
\quad \mbox{in} \quad \Gamma, \label{dbc} 
\end{eqnarray}
where $u^{\text{\emph{scat}}}\in C^2(\mathbb{R}^2\setminus 
\overline{\Omega})\cap C(\mathbb{R}^2\setminus \Omega)$.
The scattered field must satisfy
the Sommerfeld radiation condition
\[
\lim_{r\rightarrow \infty} r^{1/2}
\left( \frac{\partial u^{\text{\emph{scat}}}}{\partial r}-
ik u^{\text{\emph{scat}}}\right)=0,\quad\:r=\left\vert \xb\right\vert \, .
\]
\end{definition}

We show in the next section how to compute
$u^{\infty}(\theta)$ from an integral representation for 
$u^{\text{scat}}$.
For the moment, we simply denote this mapping, from the curve to the 
far-field pattern as 
\begin{equation}
F(\Gamma)=u^{\infty} \label{inv_operator} \, .
\end{equation}

\begin{definition}
The {\em inverse} scattering problem consists of determining 
$\Gamma$ from $u^{\infty}(\theta)$. This requires solving
the nonlinear functional equation \eqref{inv_operator}.
\end{definition}

The forward problem is both linear and well-posed, although it does
require the solution of a partial differential equation.
The inverse problem, however, is quite different.
It is both nonlinear and ill-posed, so that the development of
robust and stable solvers remains quite challenging. 
The nonlinearity leads to a non-convex 
optimization problem and the ill-posedness requires some form
of regularization.
For discussion of the question of uniqueness, we refer the reader
to \cite{Colton,kress2007}.

Following \cite{kress2007}, we note that
several approaches have been
developed to handle the nonlinearity. They can be classified as 
{iterative} methods, {decomposition} methods, and {sampling}
methods. The iterative category includes Newton's method, Landweber
iteration and the nonlinear conjugate gradient method
\cite{Colton,hohage,kress_hybrid1}. Decomposition methods proceed by first
finding an equivalent source representation that leads to the measured
far-field pattern with sources located inside $\Omega$. (The suitable placement
of sources requires some prior knowledge about the scatterer.) 
Once this linear (but ill-posed) 
problem has been solved, it remains to solve the well-posed, but nonlinear,
problem of determining $\Gamma$, using the fact that the total field
$u^{\text{inc}} + u^{\text{\emph{scat}}}$ vanishes there
\cite{kirsch_potential,Potthast_apoint-source,Potthast_apoint-source1,potthast2001point}. 
Finally, in sampling methods, one seeks to compute an ``indicator" function 
which vanishes inside $\Omega$. This category includes 
the linear sampling method of Colton and Kirsch \cite{Colton_kirsch}, 
the singular source method of Potthast 
\cite{Potthast_single,potthast2001point}, the factorization method 
of Kirsch \cite{Kirsch_2}, and the probe method of 
Ikehata \cite{Ikehata}.
We refer the interested reader to the tutorial article \cite{kress2007} and 
the text \cite{Colton} for further background material. 
The ill-posedness of the inverse problem is usually addressed by using 
some variant of Tikhonov regularization for the linearized problems 
which arise in each of the schemes listed above
\cite{cakoni2006,Colton,kaipio2010,kirsch_book}. 

In this paper, we present a Newton-like iterative method for the 
inverse obstacle scattering problem with several new features.
First, unlike many obstacle scattering
algorithms, we do {\em not} assume the boundary is star-shaped.
Instead, we assume it is bandlimited as a function of 
arclength but otherwise unconstrained. 
Second, instead of Tikhonov regularization, we simply enforce
a band-limit on the reconstructed curve using the method of
\cite{rokhlin_beylkin}. 
Third, we make use of efficient, high-order forward modeling
capabilities.
After describing the single frequency reconstruction problem, 
we investigate the quality of reconstruction using multiple angles 
of incidence as well as multiple frequencies. In the 
latter case, we rely on the recursive linearization approach developed
by Chen \cite{Chen}, 
and employed previously by Sini and Th\`{a}nh \cite{Sini11} for 
impenetrable sound-soft obstacles (the same problem considered here).
Recursive linearization was also studied by
Bao and Triki \cite{Bao} for the solution of the inverse medium problem.

In Section \ref{sec:direct},
we briefly review the direct scattering problem and its solution using
integral equation methods.
In Section \ref{sec:inv}, we describe the inverse scattering problem 
in more detail for a single frequency,
its linearization and its solution using a 
Newton-like method. We also describe our regularization based on 
band-limited approximation. 
In section \ref{sec:inv1freqmultdir}, we consider multiple frequency
data and the recursive linearization approach.
In Section \ref{sec:numerics}, we illlustrate the performance of our
method, while section \ref{sec:concl} contains some concluding remarks
and a brief discussion of future directions for research.

\section{The direct scattering problem}
\label{sec:direct}

In this section, we briefly review
the solution of the forward acoustic scattering problem for sound-soft 
impenetrable objects in two dimensions.
We refer the read to \cite{Colton} for a proof of uniqueness. We simply
note here that, in addition to being unique, the solution depends
continuously on the data $u^{\text{inc}}$ in the maximum norm.

We will make use of a boundary integral approach, since this 
requires discretization of $\Gamma$ alone and
permits the imposition of the exact Sommerfeld radiation condition.
We refer the interested reader to \cite{Colton,Nedelec} for a review 
of the relevant potential theory.

\begin{definition}
Given a boundary $\Gamma$ in $\mathbb{R}^2$, 
the single layer potential is defined by 
\begin{equation*}
\left[ S_{\Gamma,k} \varphi
\right]
(\xb):=\int_{\Gamma}\!G_k(\xb,\yb)\varphi(\yb)\;\D s (\yb),
\end{equation*}
and the double layer potential is defined by 
\begin{equation*}
\left[ D_{\Gamma,k} \varphi
\right]
(\xb):=\int_{\Gamma}\!\dfrac{\partial G_k(\xb,\yb)}{\partial \nu(\yb)} 
\varphi(\yb)\;\D s (\yb).
\end{equation*}
The normal derivative of the single layer operator is 
defined by 
\begin{equation*}
\left[ S'_{\Gamma,k} \varphi \right]
(\xb):=\int_{\Gamma}\!\dfrac{\partial G_k(\xb,\yb)}{\partial \nu(\xb)} 
\varphi(\yb)\;\D s (\yb).
\end{equation*}
\end{definition}

Here, $G_{k}(r) = \frac{i}{4} H_0(k r)$
is the Green's function for the Helmholtz equation
that satisfies the outgoing radiation
condition, where $H_0$ denotes the Hankel function of the first kind.
We note that
$S_{\Gamma,k}$ is weakly singular and the integral
is well-defined for any $\xb$. The limiting value of $D_{\Gamma,k}$ depends 
on the side of $\Gamma$ from which $\xb$ approaches the curve. 
For $\xb \in  \Gamma$,
$D_{\Gamma,k}$ is defined in the principal value sense.
Note that $S'_{\Gamma,k}$ is the adjoint of $D_{\Gamma,k}$ and
should also be interpreted in a principal value sense when 
$\xb \in \Gamma$.

\begin{rem}
When the frequency $k$ under consideration is clear from context, 
we will write $S_\Gamma, D_\Gamma, S'_\Gamma$ instead of 
$S_{\Gamma,k}, D_{\Gamma,k}, S'_{\Gamma,k}$.
\end{rem}

Using the asymptotic behavior of the Hankel function as $r \rightarrow \infty$,
we have simple formulas for the far-field patterns of the single 
and double layer potentials \cite{Colton}.
The far-field pattern of the single layer potential is given by
\begin{equation}
\left[S_{\Gamma,k}^{\infty}\varphi
\right]
(\theta):=\dfrac{e^{i\pi/4}}{\sqrt{8\pi k}}\int_{\Gamma}e^{-ik
(\cos \theta, \sin \theta) \cdot \yb}\varphi(\yb)\;\D s(\yb), \quad \theta \in U.
\label{slpfar}
\end{equation}
$U$ here is the unit circle.
The far-field pattern of the double layer potential is given by
\begin{equation}
\left[
D_{\Gamma,k}^{\infty}\varphi
\right]
(\theta):=e^{-i\pi/4}\sqrt{\dfrac{k}{8 \pi}}
\int_{\Gamma}e^{-ik
(\cos \theta, \sin \theta) 
\cdot \yb}\left(
(\cos \theta, \sin \theta) 
\cdot \nu(\yb)\right)\varphi(\yb)\;\D s(\yb).
\label{dlpfar}
\end{equation}

\subsection{The combined field integral equation} 

Two distinct methods for solving the 
forward scattering problem will be needed in our inverse scattering
algorithm.
The first is based on representing the scattered field as a linear
combination of single and double layer potentials with the same density:

\begin{equation*}
u^\text{scat}(x):= 
\left[ D_{\Gamma,k}\varphi \right](\xb) -i \eta \left[ S_{\Gamma,k}\varphi \right](\xb), 
\quad \xb \in \mathbb{R}^2\setminus\overline{\Omega},
\end{equation*}
where $\eta$ is a coupling parameter chosen to be proportional to the frequency of the incident wave \cite{kress1983minimizing}.
Imposing the boundary condition  \eqref{dbc} and using standard jump relations
for the double layer potential \cite{Colton,Nedelec}, we obtain
for $\xb\in \Gamma$, the integral equation
\begin{equation}
\left[\left(\dfrac{1}{2}I +D_{\Gamma,k}-i\eta S_{\Gamma,k}\right)\varphi\right](\xb)=-u^\text{inc}(\xb), \label{cl_int_eq}
\end{equation}
where $I$ is the identity operator. This is a well-conditioned,
resonance free Fredholm equation of the second kind \cite{Colton}.

Once $\varphi$ is known,
we may easily compute the far-field pattern 
$u^\infty(\theta)$ by using \eqref{slpfar} and \eqref{dlpfar}:
\begin{equation}
 u^\infty(\theta)=\left[\left(D_{\Gamma,k}^{\infty}-i\eta S_{\Gamma,k}^{\infty}\right)\varphi\right](\theta). \label{cl_ff_int_eq}
\end{equation}

\subsection{Green's representation} 

A second method for solving the forward problem
can be obtained by applying Green's second identity to 
$u^{\text{scat}}$ and $u^{\text{inc}}$ separately,
together with the boundary
condition \eqref{dbc}.
A straightforward calculation shows that
\begin{equation}
u(\xb) = u^{\text{inc}}(\xb) - 
\left[ S_{\Gamma,k} {\small \frac{\partial u}{\partial \nu}} \right]
(\xb),\quad x\in\mathbb{R}^2\setminus \Gamma \label{green_eq}
\end{equation}
and, therefore,
\begin{equation}
\left[ S_{\Gamma,k}{\small \frac{\partial u}{\partial \nu}} \right]
(\xb)=u^{\text{inc}}(\xb),\quad x\in\mathbb{R}^2\setminus \Gamma. \label{first_kind_equation}
\end{equation}
This is a first kind integral equation.
To improve the conditioning of this approach,
we take the normal derivative of \eqref{first_kind_equation} as
$\xb \rightarrow \Gamma$, and
add it to \eqref{green_eq} multiplied by the factor
$i \eta$. 
This results, for $\xb \in \Gamma$, in the second kind integral equation
\begin{equation}
\left[ \left(\dfrac{1}{2}I -S'_{\Gamma,k}+i\eta S_{\Gamma,k}\right)\frac{\partial u}{\partial \nu} 
\right] (\xb)=\dfrac{\partial u^\text{inc}}{\partial \nu}(\xb)
+i\eta u^\text{inc}(\xb). \label{der_int_eq}
\end{equation}

Once \eqref{der_int_eq} is solved and $\frac{\partial u}{\partial \nu}$ 
is known, it is straightforward to see from
\eqref{green_eq} that the far field pattern is
given by
\begin{equation}
u^\infty(\theta)= -\left[ S_{\Gamma,k}^{\infty} {\small \frac{\partial u}{\partial \nu}} 
\right] (\theta),
\quad \theta \in U. \label{green_ff_eq}
\end{equation}

\subsection{Numerical solution of the direct scattering problem}

It remains to discuss the discretization and solution of 
\eqref{cl_int_eq} and \eqref{der_int_eq}.
For this, we assume the boundary $\Gamma$ is parametrized by
$\gamma:[0,L]\rightarrow\mathbb{R}^2$, with 
$\gamma(0)=\gamma(L)$. We discretize the boundary using
a Nystr\"{o}m method \cite{Colton}, with $N$ equispaced
points corresponding to $t_j = (j-1)L/N$ on the boundary
for $j = 1,\dots,N$.
We assign an unknown value $\varphi(t_j)$ or 
$\frac{\partial u}{\partial \nu}(t_j)$ at each such point and enforce the
integral equation pointwise at the same points.
If the domain $\Gamma$ and the
integrands in $S_{\Gamma,k}\varphi$, $D_{\Gamma,k}\varphi$ and 
$S_{\Gamma,k}'\varphi$ were
smooth, the trapezoidal rule would yield a dense $N \times N$ matrix
whose solution would be spectrally accurate (or 
superalgebraically convergent) \cite{DR,kress_int_eq_book}. 
Because $S_{\Gamma,k}$ and the principal value for $D_{\Gamma,k}$ and 
$S_{\Gamma,k}'$ are logarithmically
singular, however, we employ the hybrid Gauss-trapezoidal rule of order 16 
due to Alpert \cite{alpert}. Without entering into
details, we simply note here that the hybrid Gauss-trapezoidal rule
replaces the diagonal band (with bandwidth sixteen) of the matrix
generated by the trapezoidal rule with special
quadrature weights in order to achieve sixteenth order accuracy.
Other high order quadrature methods could be used equally well 
\cite{helsing,kloeckner}.

For small $N$, we use direct $LU$ factorization to solve the discretized
versions of \eqref{cl_int_eq} and \eqref{der_int_eq}.
For large values of $N$, the $O(N^3)$ complexity makes the cost of 
naive $LU$ factorization prohibitive. 
Given a well-conditioned formulation, one could use fast 
multipole-accelerated iterative solution methods such as GMRES.
While this is asymptotically optimal, we are often
interested in solving scattering problems in the same geometry with 
different right-hand sides. For this,
recently developed fast direct solvers are more powerful.
These methods make use of the fact that
off-diagonal blocks of the system matrix representing 
\eqref{der_int_eq} or \eqref{cl_int_eq} are low-rank to a user-specified
precision and we refer the reader to the literature on 
hierarchical off-diagonal low-rank 
matrices (HODLR) 
\cite{siva, amirhossein2014fast, kong2011adaptive}, 
hierarchically semi-separable (HSS) or 
hierarchically block-separable (HBS) matrices 
\cite{greengard2009fast,chandrasekaran2006fast1,chandrasekaran2006fast,
ho2012fast,martinsson2009fast,martinsson2005fast}, 
and $\mathcal{H}$ and $\mathcal{H}^2$ matrices 
\cite{bebendorf2005hierarchical,borm2003hierarchical,borm2003introduction,hackbusch2001introduction}. 
Here, we will employ the
fast direct solver of Ambikasaran and Darve \cite{siva}, 
whose cost scales as $\mathcal{O}(N \log^2 N)$ for a fixed frequency,
and which performs well even for problems that are
many wavelengths in size. 

Once the densities $\{ \varphi(t_j) \}$ or 
$\{ \frac{\partial u}{\partial \nu}(t_j) \}$ are known at the 
points $\{ \gamma(t_j)\}$, we can calculate the 
far-field pattern at $M$ points $\{ \theta_\ell, \ell=1,\ldots,M \}$ 
on the unit circle using the trapezoidal rule for the far-field operators 
$S_{\Gamma,k}^{\infty}$ and $D_{\Gamma,k}^{\infty}$. 
Since the kernels of these operators are smooth (analytic), 
the trapezoidal rule yields
spectral accuracy \cite{kress_int_eq_book}. $M$ is modest in our examples
here, so we compute the far field directly, requiring 
$\mathcal{O}\left(N M\right)$ work. For larger-scale problems, fast
algorithms are available to reduce the computational complexity
\cite{candes,strain,oneil}.

\section{The inverse scattering problem}

With fast and accurate solvers available for
the forward scattering problem, we turn to the development of an 
iterative method for the solution of the inverse problem
\eqref{inv_operator}
\[
F(\Gamma)=u^{\infty} \, . 
\]

That is, given the far-field pattern, we wish to determine the shape $\Gamma$
of the scatterer itself.
We will begin with a formulation of the problem for a single incident wave at
a single frequency.

\subsection{Inversion with a single incident direction} \label{sec:inv}

Let us assume for the moment that the incident plane wave is specified by the 
direction $d$ at a single frequency $k$.
As noted in the introduction, even the theoretical foundations for the 
inverse problem are rather complicated, with issues such as uniqueness
and stability still active areas of research
\cite{Colton,kress2007}. 
Mathematically, the problem is that
the operator $F$ is a smoothing operator.
More precisely \cite{Colton}, 
$F:C^1(\Gamma)\rightarrow L^2(U)$, where $U$ is the unit circle
is continuous, compact and Fr\'echet differentiable.
Thus, inverting \eqref{inv_operator} is classically ill-posed.
Physically, the problem is that subwavelength features give rise
to evanescent waves which decay exponentially in space and are not detectable
in finite precision from $u^\infty$ (essentially, a version of 
the Heisenberg uncertainty principle).

Before turning to the question of regularization, we present an 
informal description of our iterative procedure---a 
Newton-like method to solve for the unknown $\Gamma$, 
based on the approximation
\[
F(\Gamma_j + P_j) \approx
F(\Gamma_j) + F'(\Gamma_j) P_j =u^{\infty} \, ,
\]
where $\Gamma_j$ is the $j$th guess for $\Gamma$, 
$F'(\Gamma_j)$ denotes the Fr\'echet derivative of $F$,
and $P_j$ is the update. 
The $(j+1)$st iterate is then given by 
\begin{equation}
 \Gamma_{j+1} = \Gamma_j + P_j \, . 
\label{updatecurve}
\end{equation}

\begin{rem}
We note that the Fr\'echet derivative of a compact operator is itself 
compact, so that $F'$ will inherit much of the ill-posedness of $F$.
\end{rem}

The approximation above leads to the linearized problem:
\begin{equation}
F'(\Gamma_j) P_j =u^{\infty} - F(\Gamma_j)  \, ,
\label{newton_eq}
\end{equation}
whose solution turns out to be fairly straightforward, 
due to a theorem of Kirsch \cite{Kirsch}.

\begin{theorem} \label{theo_derivative_f} 
Let $\Gamma_j$ denote a closed boundary in 
$\mathbb{R}^2$ with $u^{\text{inc}}$ a given incoming field,
parametrized by $\gamma_j(t)$ in arclength 
for $t \in [0,L_j]$.
Let $u_j^{\text{scat}}$ denote the solution to the corresponding
{\em forward} scattering problem, with 
Neumann data for the total
field on $\Gamma_j$ given by $\frac{\partial u_j}{\partial \nu}$ 
and far-field pattern $F(\Gamma_j)$.
Let $v$ denote the solution to the  
{\em forward} scattering problem with Dirichlet data
\begin{equation*}
v(t) =-\nu_j(t)\cdot P_j(t) \frac{\partial u_j}{\partial \nu}(t) 
\quad \text{on} \quad \Gamma_j,
\end{equation*}
where $P_j(t)$ is some two-dimensional perturbation of the boundary and
$\nu_j(t)$ denotes the normal vector to $\Gamma_j$. Finally,
let $v^\infty$ denote its far-field pattern.
Then
\begin{equation}
F'(\Gamma_j) P_j =v^{\infty}. \label{derivative_f}
\end{equation}
\end{theorem}

\begin{rem}
In order to make use of the preceding theorem, we need to first
obtain the normal derivative 
$\partial u_j/\partial \nu$. This is accomplished by solving 
the integral equation \eqref{der_int_eq} on $\Gamma_j$.
Given this normal derivative data, we denote by  
${\cal B}$ the diagonal operator which corresponds to multiplication by
$\partial u_j/ \partial \nu(t)$.
We then solve the Dirichlet problem for $v$ using eq. \eqref{cl_int_eq} 
and apply the far field operator from \eqref{cl_ff_int_eq}, yielding
\begin{equation}
\left(D^\infty_{\Gamma_j,k}-i\eta S^\infty_{\Gamma_j,k}\right)
\left(\dfrac{1}{2}I+D_{\Gamma_j,k}-i\eta 
S_{\Gamma_j,k} \right)^{-1} {\cal B} (\nu_j \cdot P_{j}) =
u^\infty -F(\Gamma_j). 
\label{full_newton}
\end{equation}
\end{rem}

Rather than letting $P_j$ be an arbitrary perturbation of the curve
$\Gamma_j$, however, we will assume that it
lies in the normal direction: 
\[ P_j(t) = \nu_j(t) p_j(t) \, , \]
where $p_j(t)$ is a scalar function. This avoids certain kinds of 
nonuniqueness
which tangential motion of the boundary would permit, making the recovery
even more ill-conditioned. Thus, eq. \eqref{full_newton} is 
replaced by
\begin{equation}
\left[ \left(D^\infty_{\Gamma_j,k}-i\eta S^\infty_{\Gamma_j,k}\right)
\left(\dfrac{1}{2}I+D_{\Gamma_j,k}-i\eta 
S_{\Gamma_j,k} \right)^{-1} {\cal B} \right]  p_{j}(t) =u^\infty -F(\Gamma_j). 
\label{fin_ful_newton}
\end{equation}
Hereafter, we will refer to the operator on the left-hand side of 
eq. \eqref{fin_ful_newton} acting on $p_j(t)$ as $F'(\Gamma_j)$.
Several important features of the iteration remain to be discussed.

\begin{enumerate}
\item We must ensure that the curve $\Gamma_{j+1}$ constructed
via \eqref{updatecurve} is not self-intersecting.
Therefore, we define
\[ \Gamma^l_{j+1}=\Gamma_{j}+\rho \lambda^l P_{j}, \] 
where $\rho$ is a user-specified parameter, and
$\lambda < 1$ provides additional damping of the Newton step if needed.
We begin with $l=0$ and accept the curve
based on oversampling the curve $\Gamma^l_{j+1}$
at a large number of points $N_s$. We check that the polygonal 
approximation of the curve is not self-intersecting. This is trivial
to do using a naive $\mathcal{O}(N_s^2)$ algorithm. More sophisticated
algorithms achieve a computational complexity of 
$\mathcal{O}(N_s\log N_s)$ \cite{Bentley1979,Preparata1985,Shamos1976}.
If $\Gamma^l_{j+1}$ fails to be simple, we increase $l$ and repeat.

\item
To overcome the ill-posedness of the problem, we will restrict $p_j$
to be a band-limited curve:
\begin{equation}
p_j(t) = \sum_{m=-b/2}^{b/2-1} p_{j,m} e^{i m 2 \pi t/L} \, . 
\label{pdiscrete}
\end{equation}
We will return to the selection of the parameter
$b$ in section \ref{bandlim}.

The integral equation \eqref{fin_ful_newton} then takes the form
of a ``discrete-to-continuous" map, with $b$ degrees of freedom
which we can solve in a least squares sense.
The fully discrete version of \eqref{fin_ful_newton} is:
\[ \left[ A^\infty C^{-1} B O \right] \vecp = R \,  \]
where $A^\infty$ is an $M\times N$ discretization of 
$\left(D^\infty_{\Gamma_j,k}-i\eta S^\infty_{\Gamma_j,k}\right)$,
$C$ is an $N\ \times N$ discretization of 
$\left(\dfrac{1}{2}I+D_{\Gamma_j,k}-i\eta S_{\Gamma_j,k} \right)$,
$B$ is a diagonal $N \times N$ matrix discretizing ${\cal B}$,
$O$ is an $N\times l$ matrix mapping the coefficients $\vecp$ in 
\eqref{pdiscrete} to equispaced samples of $p_j(t)$ on $\Gamma_j$
and $R$ is a discretization of the right-hand side
$u^\infty -F(\Gamma_j)$ at $M$ points.
(We discretize the integral eq. operator $C$ from \eqref{cl_int_eq} 
as described in section \ref{sec:direct}.)
Thus, in the fully discrete version, we have
$F' =  A^\infty C^{-1}\, B\, O$, which  
is of dimension $M \times b$.
We assume that $M$ and $b$ are sufficiently small that standard
linear algebra tools for least squares problems can be applied, 
such as the $QR$ factorization, so long as $M > b$.

\item
It is 
more efficient to compute the entries of $(F')^T =
O^T B^T C^{-T} (A^\infty)^T$ than $F'$ itself, since we then need to apply
$C^{-T}$ to $M$ distinct vectors, namely the columns of $(A^\infty)^T$.
(We are assuming here that $M$ is significantly smaller than $N$.)
With the fast direct solver
of \cite{siva}, this requires $O(MN \log^2 N)$ work.
The operator $O^T$ can be applied using the fast Fourier transform
at a net cost of $O(M(b+N) \log N)$ work.

\item
After each iteration, we have a new boundary $\Gamma_{j+1}$ sampled at
$N$ points: 
$\gamma_{j+1}(t_n) = \gamma_{j}(t_n) + \rho \lambda^m P_{j}(t_n) \nu_j(t_n)$.
where $m$ is determined by the acceptance criterion discussed above.
To further improve the conditioning of the inverse problem,
we use the algorithm of
Beylkin and Rokhlin \cite{rokhlin_beylkin} to filter the 
curve $\Gamma_{j+1}$ and return
$N$ points equispaced in arclength on this resampled curve.
Without entering into details,
we have implemented the algorithm of \cite{rokhlin_beylkin} to
permit filtering of the curve beginning at frequency $b$,
with a smooth roll-off until frequency $b+N_b$. 
On output, the algorithm returns 
$N$ points equispaced in arclength with respect to the new (smooth) 
parameterization.

\item
In our examples below, we choose $\Gamma_0$ to be the unit disk.
Better initial guesses can be obtained by using 
a sampling method \cite{Colton_kirsch,Ikehata,Kirsch_2,Potthast_single,potthast2001point}
which provides a good approximations to the convex hull of 
the unknown object.

\item
In addition to the use of damping to ensure that the curve is simple
(not self-intersecting), we need some stopping criteria for the Newton
iteration. For this, we set a maximum number of Newton steps, a minimum
norm for the solution to the linearized problem $|p_j(t)|$ and a 
residual tolerance ($\| F(\Gamma)-u^{\infty} \| < \epsilon$).

The choice of the parameter $\rho$ in the damped Newton step can also be 
important \cite{Conn,Schnabel}.
(Small values tend to require more iterations but can improve convergence.)
In \cite{li_bao} Li and Bao suggest a frequency-dependent
damping parameter. We have not explored this choice systematically
in the present paper.
\end{enumerate}

\subsection{Inversion with multiple incident directions}
\label{sec:inv1freqmultdir}

Suppose now that we seek to solve the inverse problem 
when scattering data is available from several angles of incidence. 
Let $L$ denote the number of incident plane waves:
\[ 
u_m^{\text{inc}}(\xb)=e^{ik \xb \cdot \db_l} 
\] 
with directions $\db_l$ for $l=1,\ldots,L$.
To solve this problem, 
eq. \eqref{newton_eq} is replaced at each iteration
by the system of equations:
\begin{equation}
\left[ \begin{array}{c}
F'_1(\Gamma_j) \\
F'_2(\Gamma_j) \\
\vdots \\
F'_{L}(\Gamma_j) \end{array} \right] P_j=\left[  \begin{array}{c}
u^{\infty}-F_1(\Gamma_j)\\
u^{\infty}-F_2(\Gamma_j)\\
\vdots \\
u^{\infty}-F_{L}(\Gamma_j)\end{array} \right]\label{multid_sfreq}
\end{equation}
where $F_l$ is the far-field operator on $\Gamma_j$ for
data generated by the incident
plane wave $u_l^{\text{inc}}(\xb)$ and 
$F_l'$ is its Fr\'echet derivative.

The discretization of each of the operators $F_l$ and $F'_l$
is carried out as in the previous section.
Note that,
since the unknown in the fully discrete version of \eqref{multid_sfreq}
is still $\vecp$, we can only improve the conditioning of the least squares
problem when compared to the setting with a single angle of incidence.
This makes physical sense;  illuminating the unknown obstacle from
more directions should certainly make the recovery problem easier.

\section{Multiple frequency measurements and recursive linearization}

Since the number of nontrivial measurements that 
can be made in the far field is proportional to the size of the object in 
wavelengths, it is reasonable to expect that greater resolution 
should be obtained as the frequency increases. There are fundamental
difficulties, however, with using Newton's method
(or one of its variants) for inverse scattering at a single frequency $k$.
Mainly, when $k$ is large, the initial guess for the boundary shape $\Gamma$
must be close to the correct solution in order for the iteration not to 
diverge or be trapped in a spurious local minimum.
Moreover, only the ``illuminated" portion of the boundary can be recovered
with fidelity, as discussed in \cite{carlos_thesis}.
When $k$ is small, however, despite the fact that the 
inverse scattering problem is ill-posed, if one only seeks to recover 
a small number of parameters (say the centroid and area of the scatterer),
a few Newton iterations are sufficient for convergence {\em without} a 
good initial guess. 

This interplay between easy recovery at low frequencies of a blurry
reconstruction and the need for a good initial guess at high frequency
to achieve higher fidelity reconstruction led Chen \cite{Chen}
to introduce the {\em recursive linearization algorithm} (RLA).
The essential idea is that given the (blurry but converged) 
solution to the inverse scattering problem at frequency $k$, one
step of Newton's method is sufficient to get the converged solution for 
frequency $k+ \delta k$, for sufficiently small $\delta k$.
That is the sense in which the problem has been linearized.
It has been used previously by Sini and Thanh \cite{Sini11} and Borges
\cite{carlos_thesis} for the 
inverse obstacle problem for a single angle of incidence with
an unknown star-shaped object. 
The procedure is initiated at a sufficiently small frequency $k$ so that
Newton's method converges without a good initial guess.
We refer the reader to \cite{Chen} for a more detailed explanation.

When a finite number of frequency measurements are available,
$k_1 < k_2 < \dots < k_J$, a modification
of RLA is to begin at the lowest frequency $k_1$ 
and to solve the inverse scattering
problem to obtain a guess which we will denote by $\Gamma^{(k_1)}$.
$\Gamma^{(k_1)}$ can then be used as the initial guess for a new Newton
iteration at frequency $k_2$, with the recursion continuing until the 
maximum frequency $k_J$ has been reached.
A formal theory, of course, might require tight control on the spacing in
frequency $\delta k = k_{j+1} - k_j$, and we use the RLA framework here
simply as a sensible guideline, as in 
\cite{carlos_thesis,Sini11}. 

In summary, for the muti-frequency, multi-direction problem,
we assume we are given
measurements of the far-field pattern at $M$ angular points $\theta_l$, 
generated by the scattering 
of $J \cdot L$ incident plane waves 
$u^{\text{inc}}_{j,l}(x)=e^{ik_j \xb\cdot \db_l}$ 
for $l = 1,2,\ldots,L$ and $j=1,2,\ldots,J$,
where $k_1<k_2<\cdots<k_{J}$. The goal, as always, is to reconstruct
the unknown object $\Omega$ with boundary $\Gamma$.

\subsection{Setting the bandlimit} \label{bandlim}

In order for the inverse problem
to be well-conditioned at each frequency, it is 
important to choose the parameter $b$ in \eqref{pdiscrete} correctly.
The value of $b$ determines how complicated a curve we are able to 
reconstruct since it controls the complexity of the perturbation $P_j(t)$.
Clearly, its value should be approximately 
proportional to $k$ in order to achieve the best possible resolution.
Smaller values would lead to stable but excessively pessimistic 
reconstructions while larger values would permit subwavelength
instabilities to develop.

\section{Numerical experiments} \label{sec:numerics}
In this section, we illustrate the performance
of our inverse obstacle scattering algorithm.
Before turning to specific examples, a few general remarks are in order.
First, with a single incident wave, we should not expect to obtain a 
good reconstruction of the poorly illuminated ``backside" of the boundary.
(This effect is more pronounced at higher 
frequencies \cite{carlos_thesis,Sini11}.) 
With multiple directions of incidence, however, it should be possible to 
reconstruct the entire object with much greater fidelity.
Second, because we are carrying out our reconstruction over the space
of band-limited closed curves, geometric singularities and regions of 
high curvature will inevitably be filtered (with the extent of filtering
a decreasing function of frequency.)

For each object, we construct simulated measurement data in the far field 
pattern by solving the direct scattering problem, as discussed in 
Section \ref{sec:direct}. We compute the solution
to \eqref{der_int_eq} or \eqref{cl_int_eq} with at least ten digits of 
accuracy to generate far field data
$u_{k,\db}^{\infty}$ and then introduce Gaussian noise in the data. 
Here $k$ denotes the frequency and $\db$ the direction of incidence.
We set
\begin{equation*}
v^{k,\db}_{\infty} = u_{k,\db}^\infty+
\delta \dfrac{\|u_{k,\db}^\infty\|}{\|\epsilon_1 + i \epsilon_2\|}
(\epsilon_1 + i \epsilon_2),
\end{equation*}
where $\epsilon_1$ and $\epsilon_2$ are normally distributed variables 
with mean zero and variance one.
As in \cite{Ivanyshyn}, we also make sure to use a different number
of quadrature points in the discretization of the  
the direct problem and in the integral equations used in solving the 
inverse problem.

In our examples, we use $\delta=0.05$, 
and let
$\db_l=2\pi l/L$ (for $l=1,\ldots,L$) for wavenumbers 
$k_j=k_0+j \Delta k$ with $j=1,\ldots,J$.
The far field pattern is obtained at the angles 
$\theta_\ell=(2\ell-1)\pi/M$, for 
$\ell=1,\ldots,M$. 
At each wavenumber $k$, the number of quadrature points used 
in solving the direct scattering problem is $N=\ceil{100k|\Gamma|}$, 
where $|\Gamma|$ is the perimeter of the object. 
The number of collocation points used in solving 
the inverse problem is denoted by $N_1$.
The specific values for $M$, $L$, $J$, $k_0$, $\Delta k$ and $N_1$ 
are described for each numerical experiment below. 

\begin{example}{A star-shaped object with 7 oscillations}
\end{example}

\begin{figure}
  \centering
\subfigure[A single incident wave]{
      \includegraphics[scale=.3]{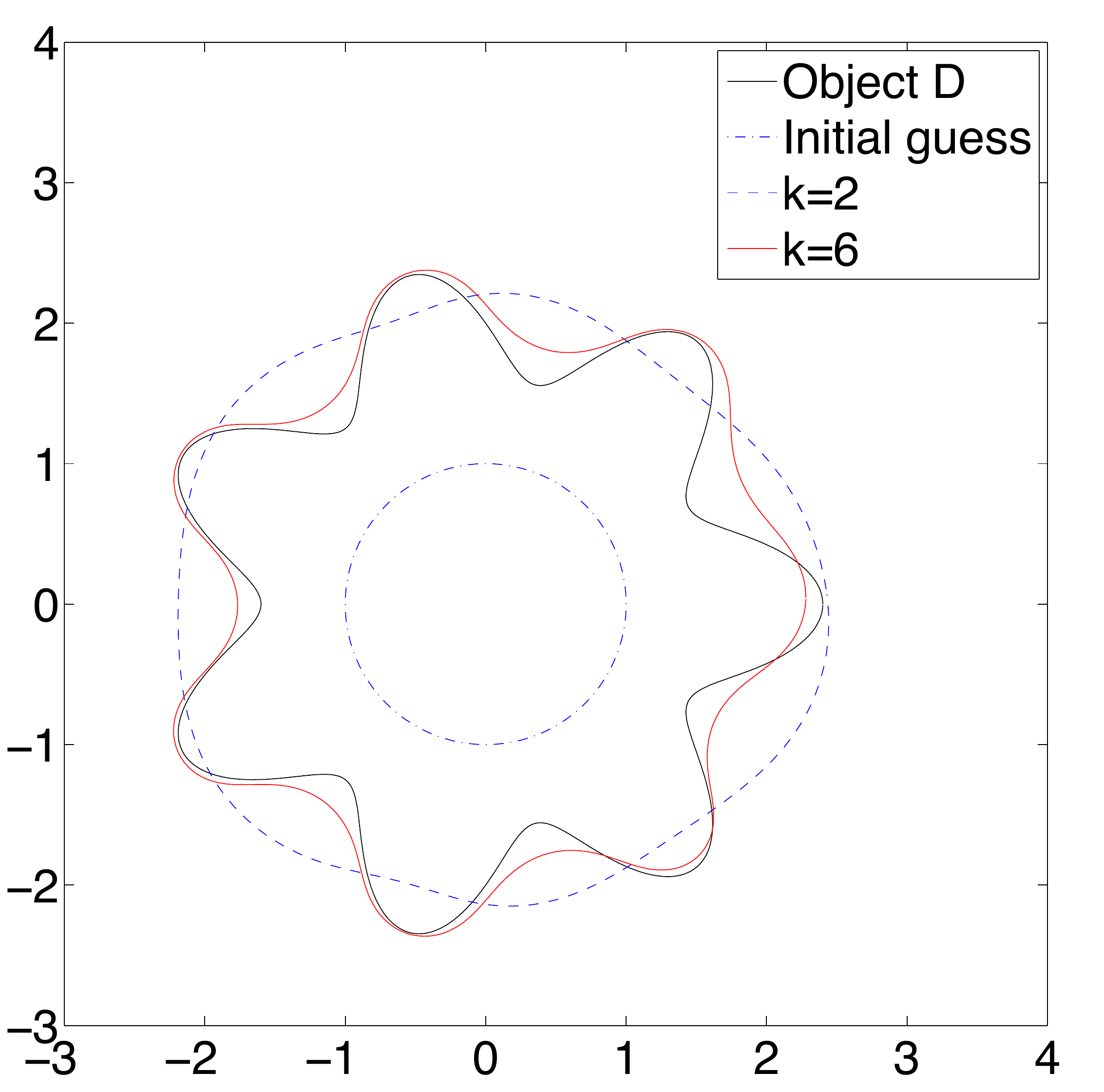}\label{gear1}
      }
  \subfigure[Four incident waves]{      
            \includegraphics[scale=.3]{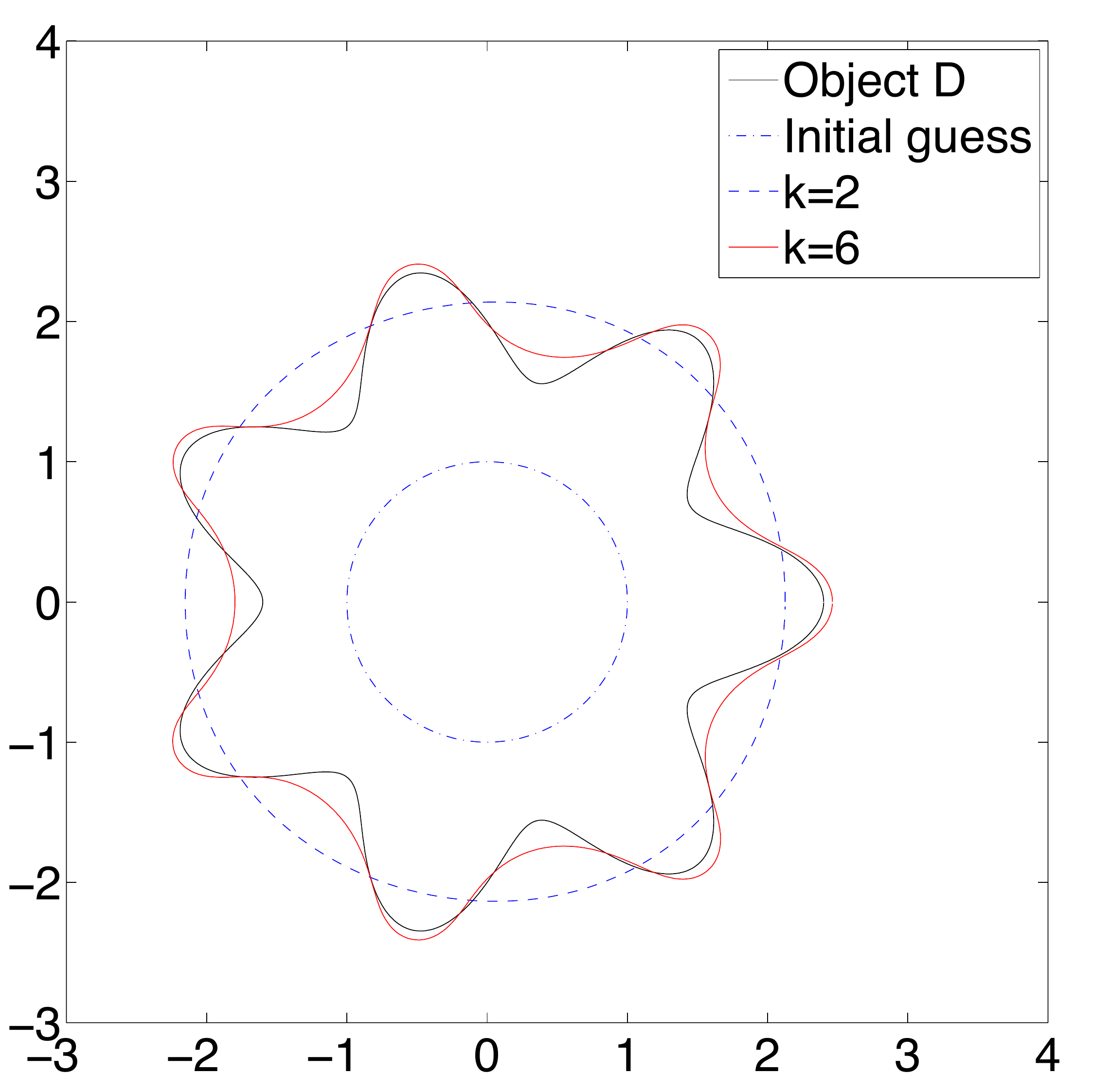}\label{gear4}
  }
  \caption[Applying the method to reconstruct the star-shaped object of Example 1.]{\footnotesize 
Reconstruction using recursive linearization with Newton's method for 
the object in Example 1 for $L=1$ (a) and for $L=4$ (b).} \label{nnss}
\end{figure}

In our first example, we consider a star-shaped object with 
parameterization $\gamma:[0,2\pi]\rightarrow\mathbb{R}^2$ 
given by
\begin{equation}
\gamma(t)=\left(2+0.2\cos(7t)\right)\left(\cos(t),\sin(t)\right).
\end{equation}
The initial wavenumber is $k_0=0.5$, $\Delta k=0.5$, and $J=11$. 
We let $M=32$ and use the damping parameter 
$\rho=0.1$ for wavenumbers up to $k=5$ and $\rho=0.1/k$ for higher wavenumbers. 
Within Newton's method, the stopping criterion used 
for the nonlinear residual is
$\epsilon = \| F(\Gamma) - u^\infty \| < 10^{-4}$.
We used $N_1=\ceil{10k|\Gamma_j|}$ quadrature points in solving
the integral equations within the Newton step. 
In Fig. \ref{gear1}, we plot the reconstruction using 
a single direction of incidence, $L=1$, and use the bandlimit $b=\ceil{k}$ 
at each wavenumber, while in Fig. \ref{gear4} 
we plot the reconstruction for $L=4$ and use the bandlimit $b=2\ceil{k}+1$. 
In both cases, the filter in the curve resampling algorithm
uses a smooth roll-off from $b$ to $b+50$. In both 
cases, we show the solutions obtained for $k=2$ and $k=6$.
The reason we can choose a larger value for
$b$ when using multiple angles of incidence is that the system matrix
\eqref{multid_sfreq} remains well-conditioned, as noted in section 
\ref{sec:inverse}.

To illustrate the advantage of solving the forward scattering problem 
using the integral equation \eqref{cl_int_eq}
by applying the HODLR algorithm, we compare the run-time 
of HODLR and direct factorization for problems with increasing wavenumber $k$ 
and number of discretization points, using the
star-shaped object with 7 oscillations as the geometry. The error is computed
by taking as boundary data the field induced by a singular point source at the origin,
which lies in the interior of the scatterer and testing the computed exterior solution
at the target point $(10,8)$.
 
\begin{table}
\caption{Running time for direct scattering problem}\label{forward_tab}
\begin{center}
\begin{tabular}{c|c|c|c|c|c}
\hline
\multicolumn{1}{ c| }{\multirow{2}{*}{\hspace{0.5cm}$k$\hspace{0.5cm}} } &
\multicolumn{1}{ |c| }{\multirow{2}{*}{\hspace{0.5cm}$N$\hspace{0.5cm}} } & \multicolumn{2}{ |c| }{Direct solver} & \multicolumn{2}{ |c }{HODLR} \\ \cline{3-6}
\multicolumn{1}{ c|  }{}                        & \multicolumn{1}{ |c|  }{} & Time(s)  & Error & Time(s) & Error    \\
\hline\hline
1     & 360 & 0.17 & 7.711837e-13 & 0.72 & 7.715468e-13 \\
2     & 720 & 0.57 & 1.764748e-12 & 1.64 & 1.766015e-12 \\
4     & 1440 & 2.08 & 3.097516e-12 & 3.80 & 3.094913e-12 \\
8     & 2880 & 7.19 & 1.988648e-12 & 7.91 & 2.120391e-12 \\
16   & 5760 & 33.48 & 2.298391e-11 & 17.20 & 2.303258e-11 \\
32   & 11520 & 314.83 & 1.957544e-11 & 37.91 & 1.939017e-11 \\
64   & 23040 & -- & -- & 95.44 & 3.924397e-11 \\
128 & 46080 & -- & -- & 278.49 & 9.649962e-11 \\
\hline
\end{tabular}
\end{center}
\end{table}

Note that the direct solver is faster than the HODLR for small problems,
but that it scales much better with $N$. The run-times for direct factorization are omitted for 
the largest problems, since prohibitive amounts of memory would be needed.

\begin{example}{An aircraft-like object}\label{results_3}\end{example}

\begin{figure}[h!]
  \centering
\subfigure[Reconstruction with 1 incident direction]{
      \includegraphics[scale=0.28]{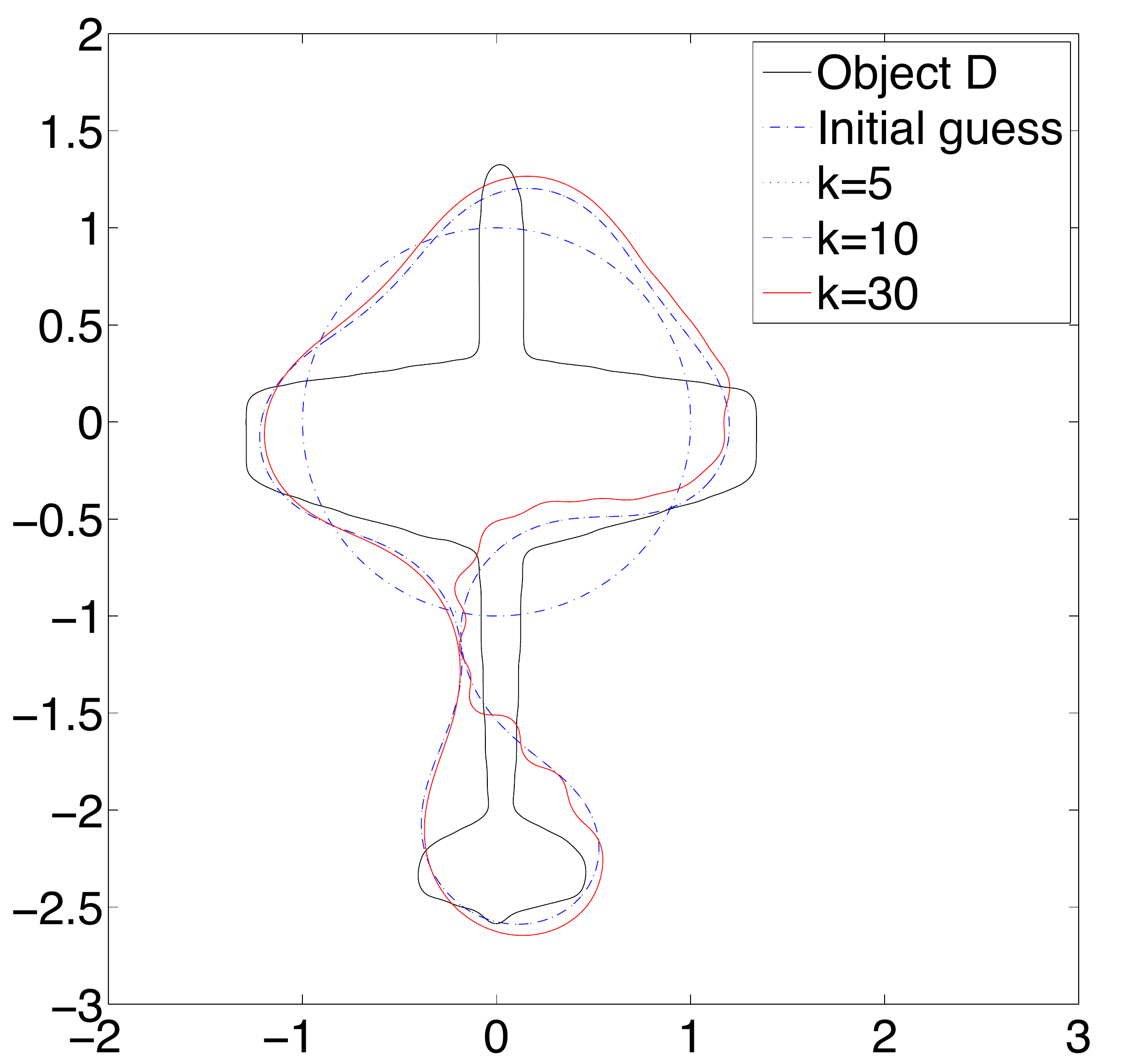}\label{planedir1}
      }
  \subfigure[Reconstruction with 2 incident directions]{      
            \includegraphics[scale=0.28]{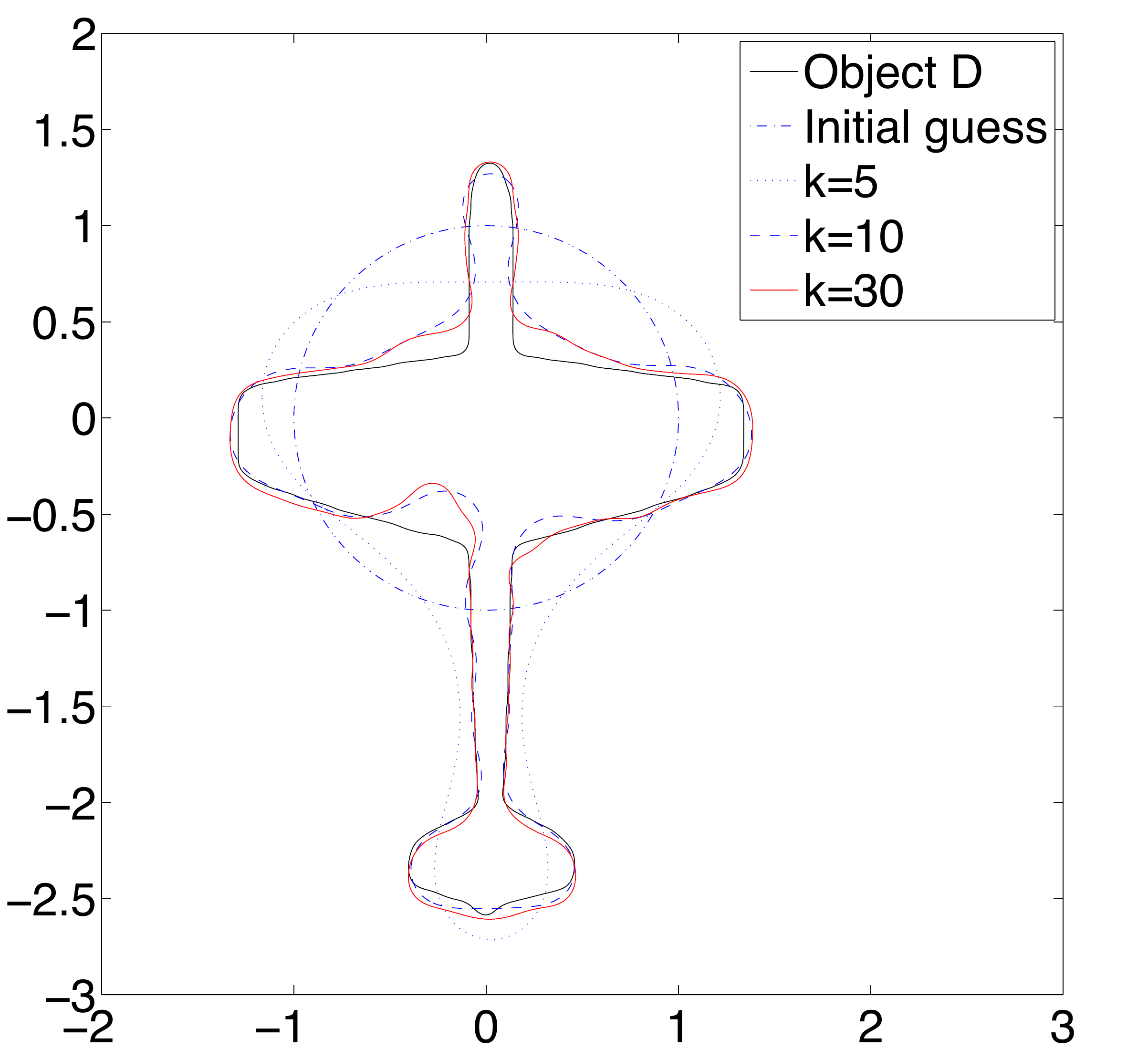}\label{planedir2}
  }

  \subfigure[Reconstruction with 3 incident directions]{      
            \includegraphics[scale=0.28]{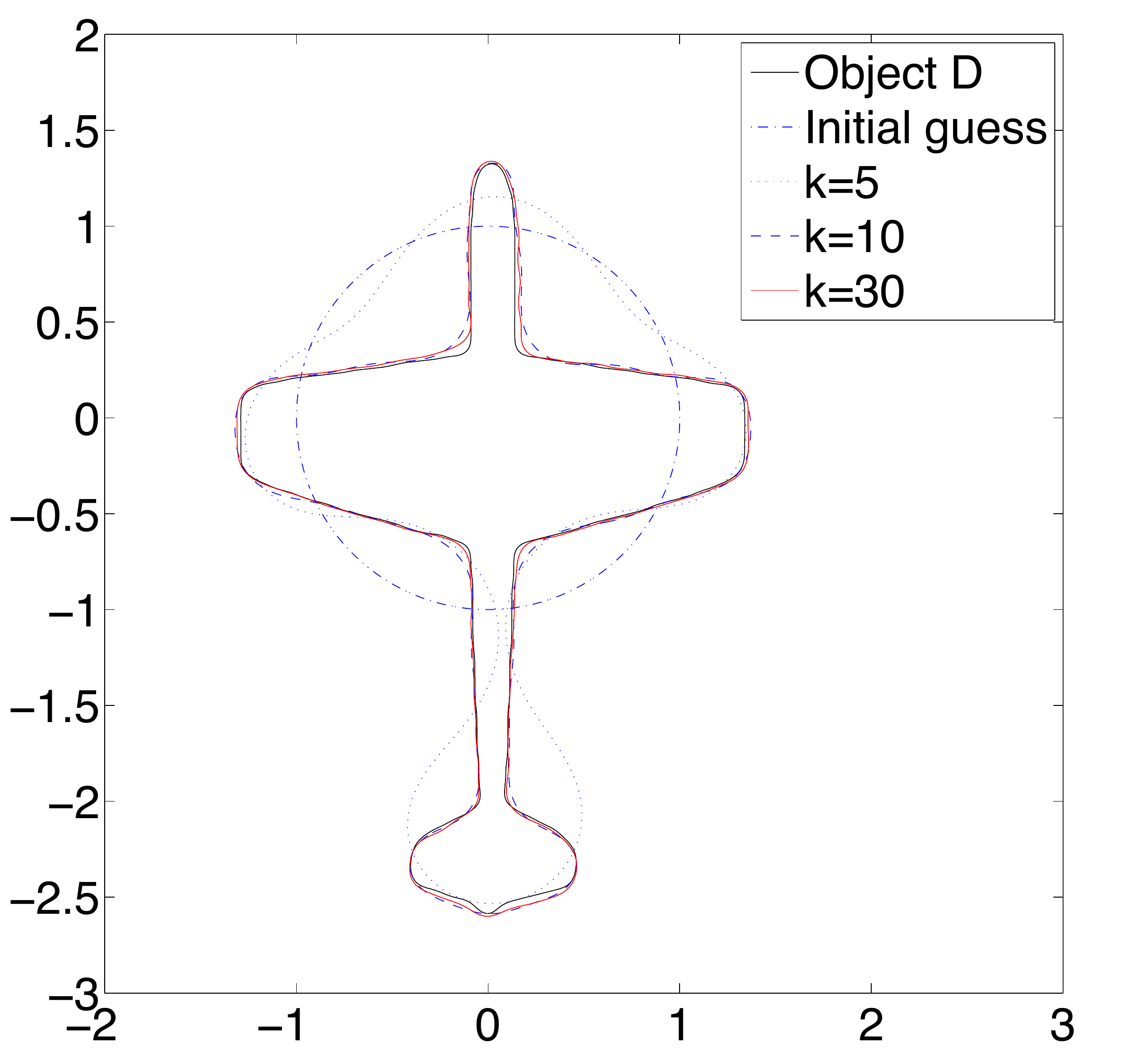}\label{planedir3}
  }
  \subfigure[Reconstruction with 4 incident directions]{      
            \includegraphics[scale=0.28]{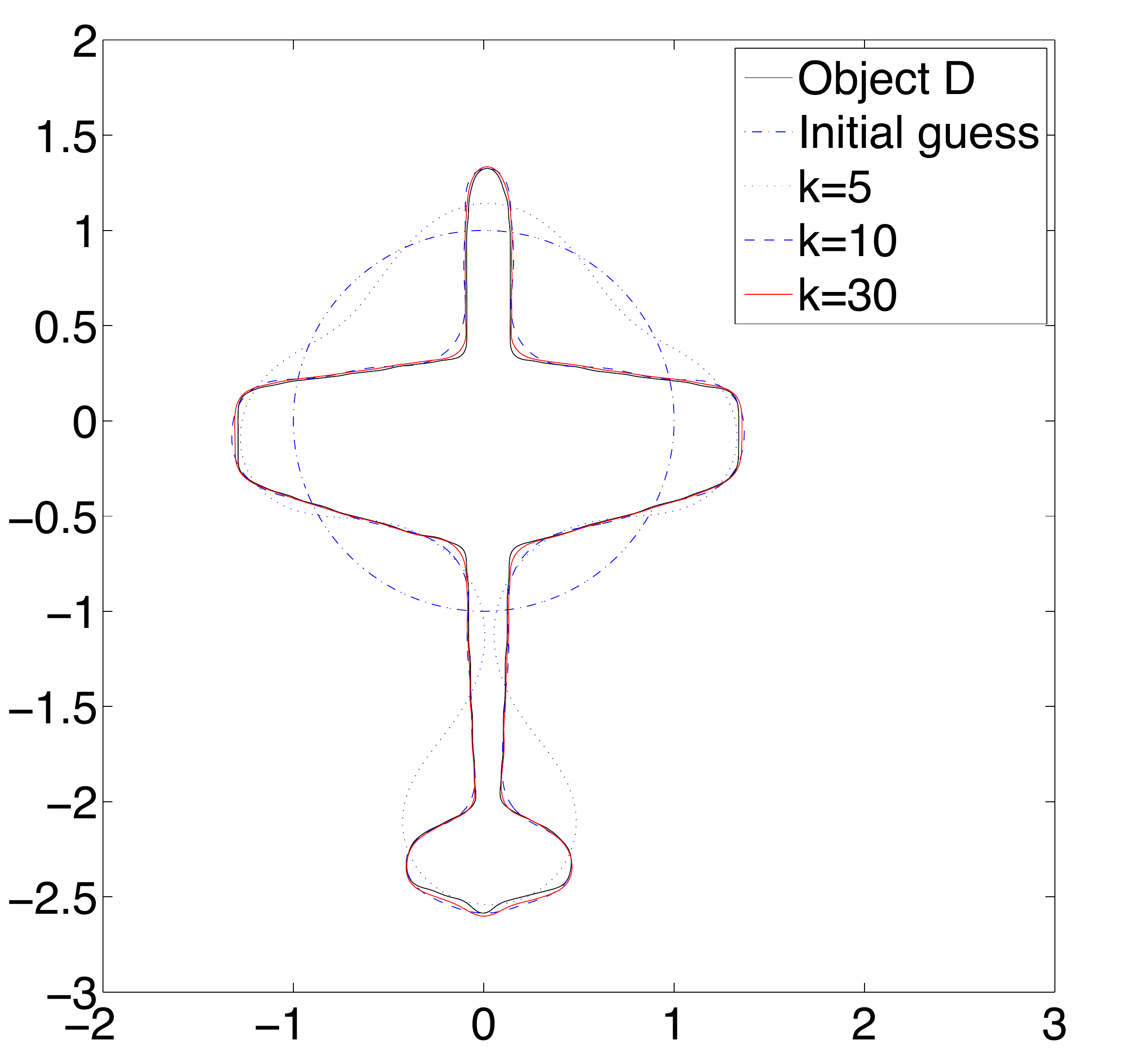}\label{planedir4}
  }

  \caption[Applying the method to reconstruct the aircraft shaped object using an increasing number
of incident directions.]{\footnotesize 
Reconstruction of the aircraft-like object using Newton's method with an increasing number of incident directions. } \label{planedirs}
\end{figure}

We next consider a more complicated (and not star-shaped) object, shown in
Figs. \ref{planedir1}-\ref{planedir4}. Our initial guess is the circle of 
radius 1 centered at the origin.
We used an initial wavenumber of $k_0=0.5$, with $\Delta k=0.5$, and $J=60$,
and let the number of far-field measurements be $M=128$. Within Newton's method, the stopping criterion used 
for the nonlinear residual is
$\epsilon = \| F(\Gamma) - u^\infty \| < 10^{-3}$.
We used $N_1=\ceil{20k|\Gamma_j|}$ quadrature points in solving
the integral equations within the Newton step and let the damping parameter $\rho=1$. 

In order to show the improvement of the reconstruction using an increasing number of incident 
plane waves, we present the final reconstructions for the cases $L=1,2,3,4$ with $b=2\ceil{k} +1$
(Figs. \ref{planedir1} - \ref{planedir4}).

\begin{figure}[h!]
  \centering
\subfigure[Reconstruction of the aircraft]{
      \includegraphics[scale=.3]{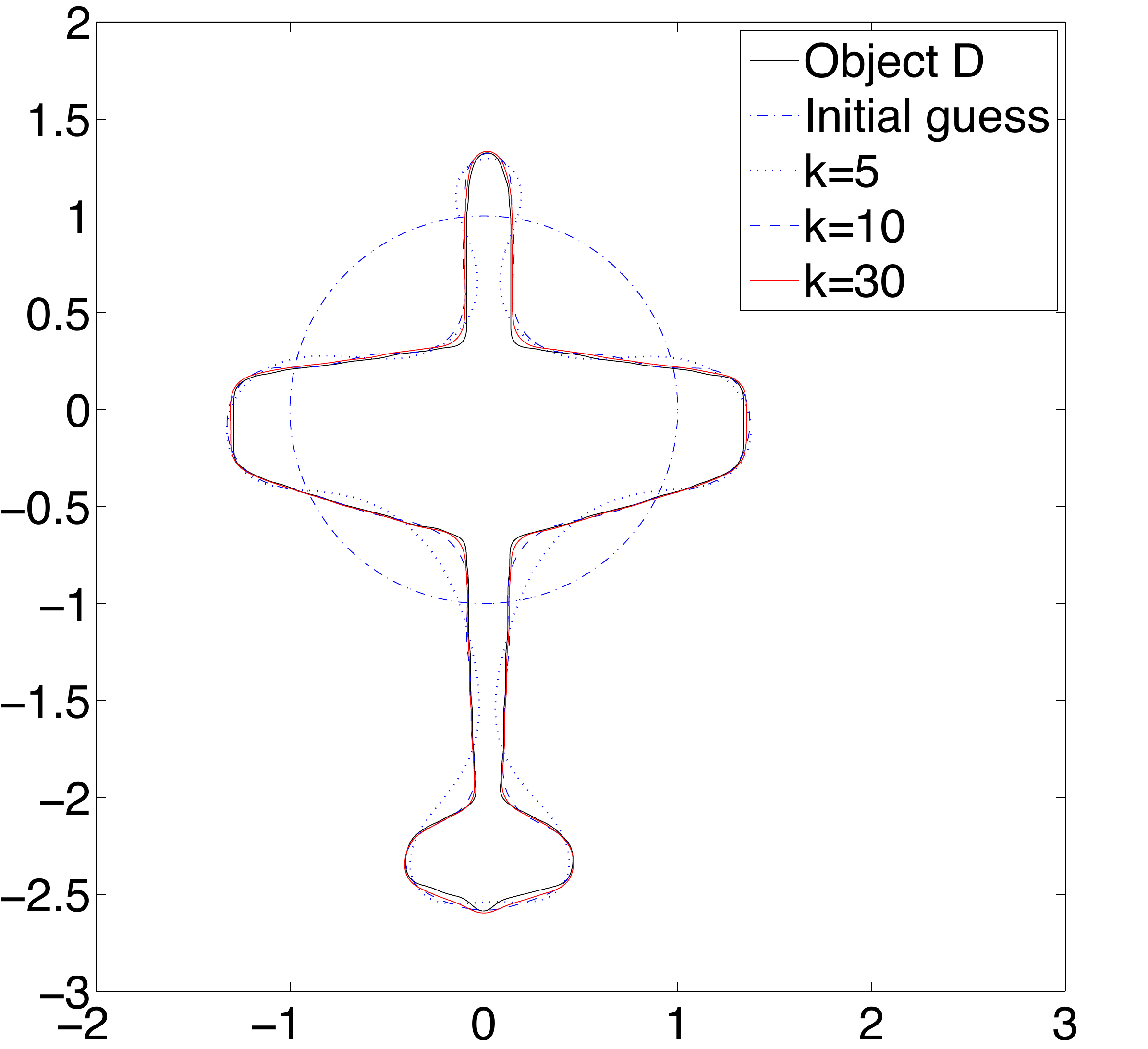}\label{plane1a}
      }
  \subfigure[Top corner of the intersection of the right wing and the fuselage of the aircraft]{      
            \includegraphics[scale=.3]{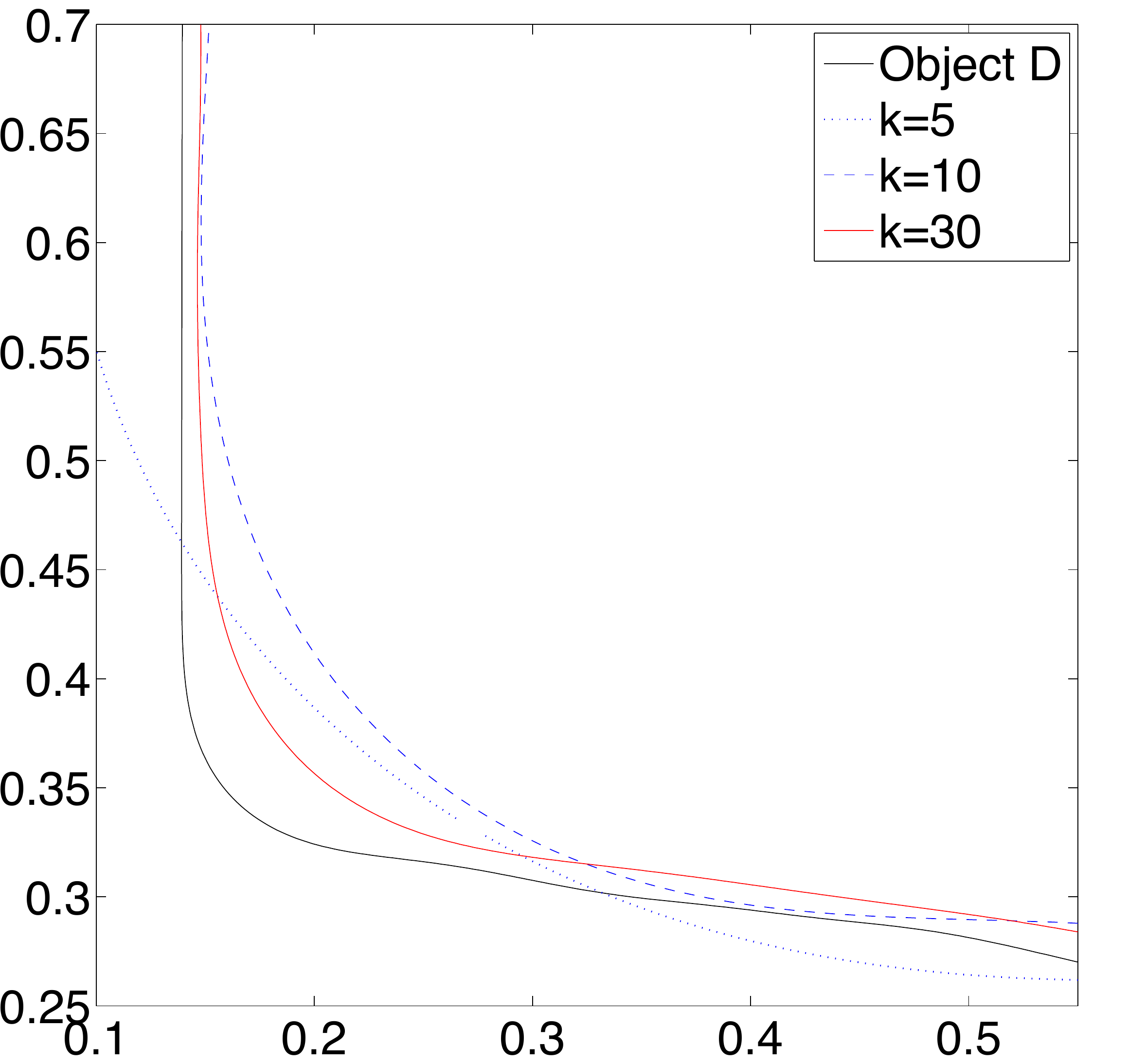}\label{plane1b}
  }
  \subfigure[Bottom corner of the intersection of the right wing and the fuselage of the aircraft]{
      \includegraphics[scale=.28]{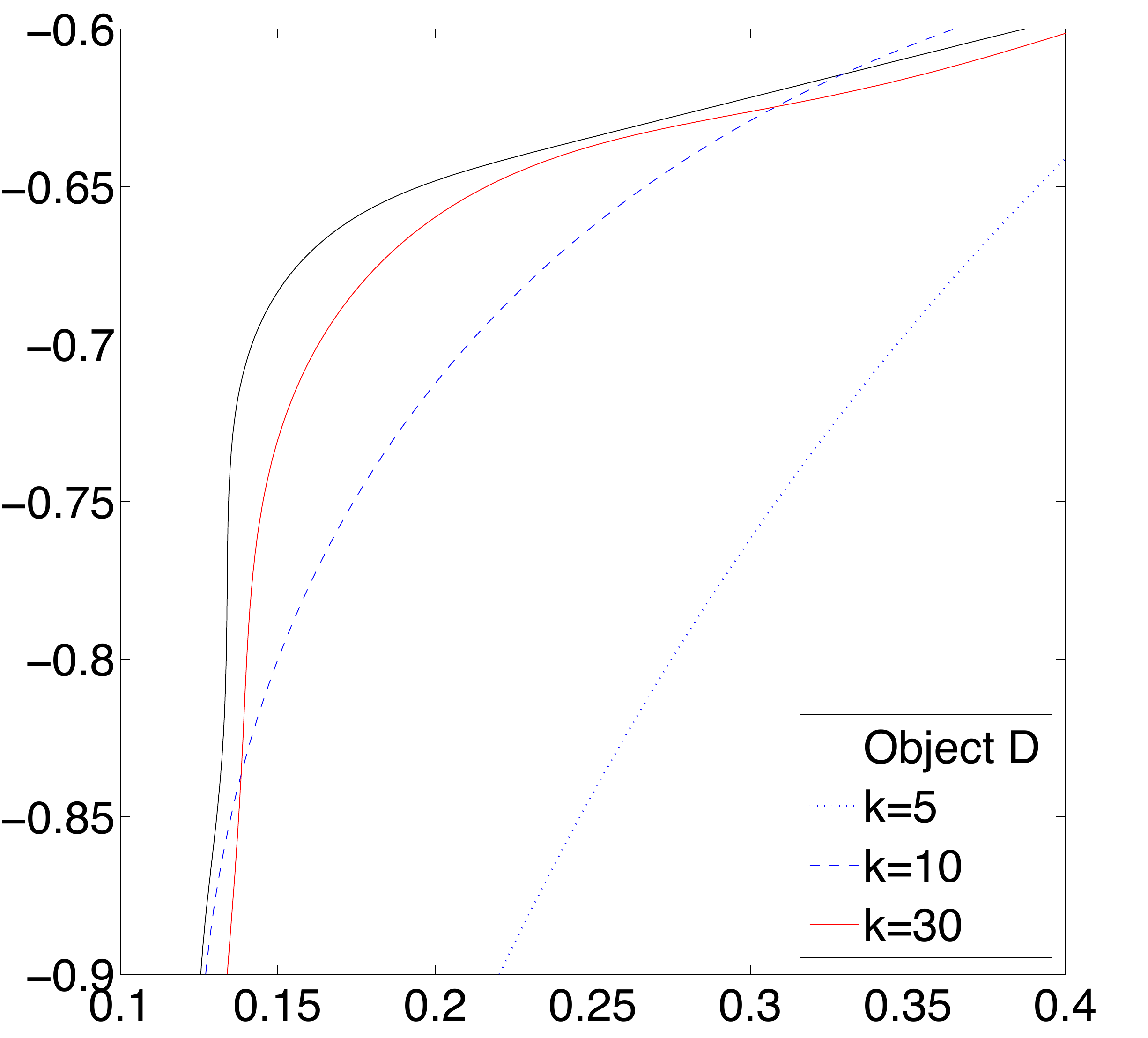}\label{plane1c}
      }
  \subfigure[Top right corner of the tail of the aircraft]{      
            \includegraphics[scale=.28]{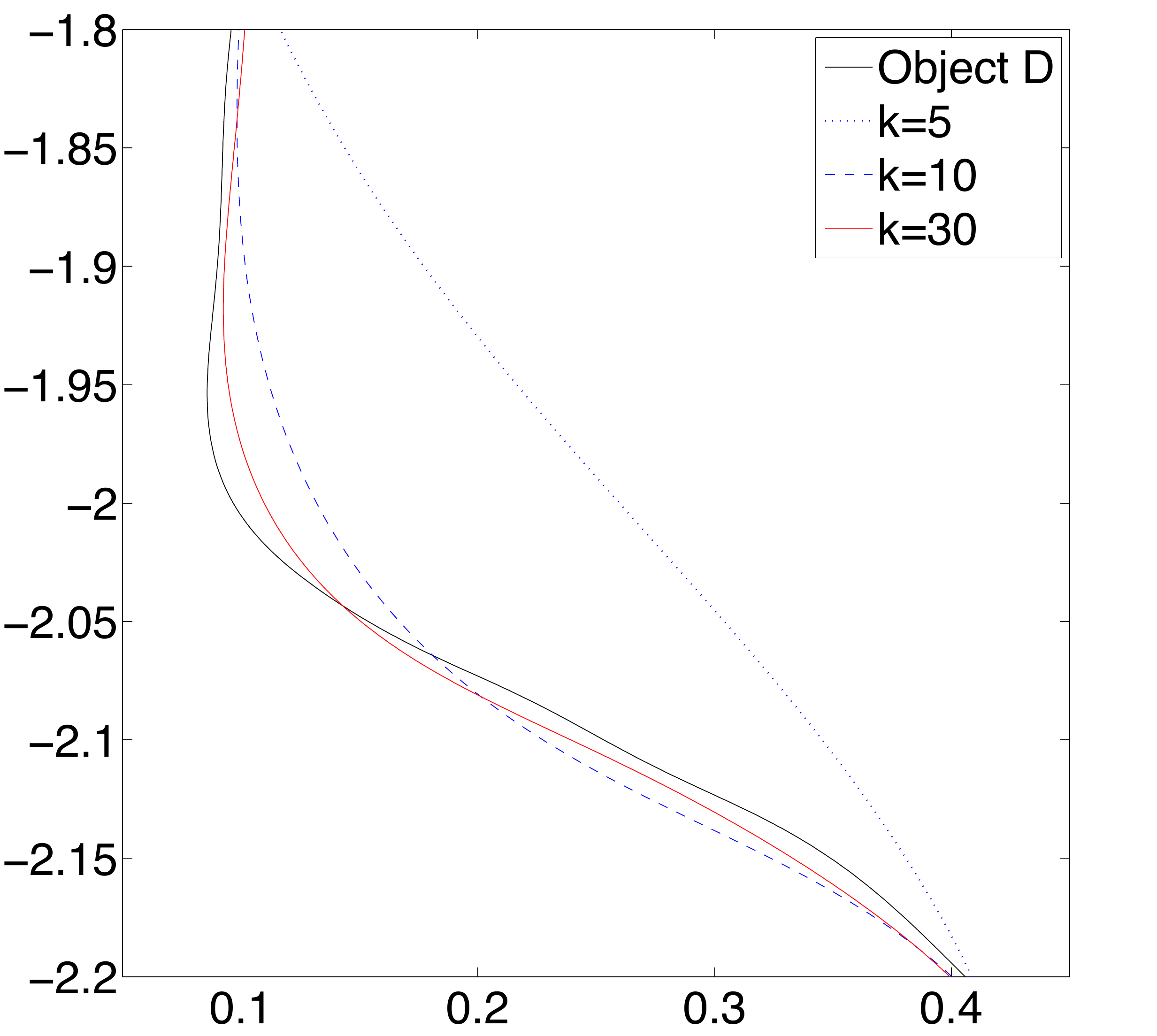}\label{plane1d}
  }

  \caption[Applying the method to reconstruct the aircraft shaped object.]{\footnotesize 
Reconstruction using recursive linearization with Newton's method to 
reconstruct the aircraft-like object of Example 2.} \label{plane1}
\end{figure}

Finally, we present a reconstruction using six directions of incidence and the bandlimit $b=2\ceil{k}+ 1$, 
with a smooth roll-off in the curve resampling from $b$ to $b+130$. All other parameters are chosen as before. 
In Fig. \ref{plane1a}, 
we plot the solution obtained with our scheme for varying maximum 
frequencies. In Figs.  \ref{plane1b}--\ref{plane1d}, we present zoomed in details
of the reconstruction in certain high-curvature regions 
of the object.

\begin{example}{Reconstruction of an elongated object}\label{results_4}\end{example}

\begin{figure}[h!]
  \centering
\subfigure[Reconstruction of the submarine]{
      \includegraphics[scale=.3]{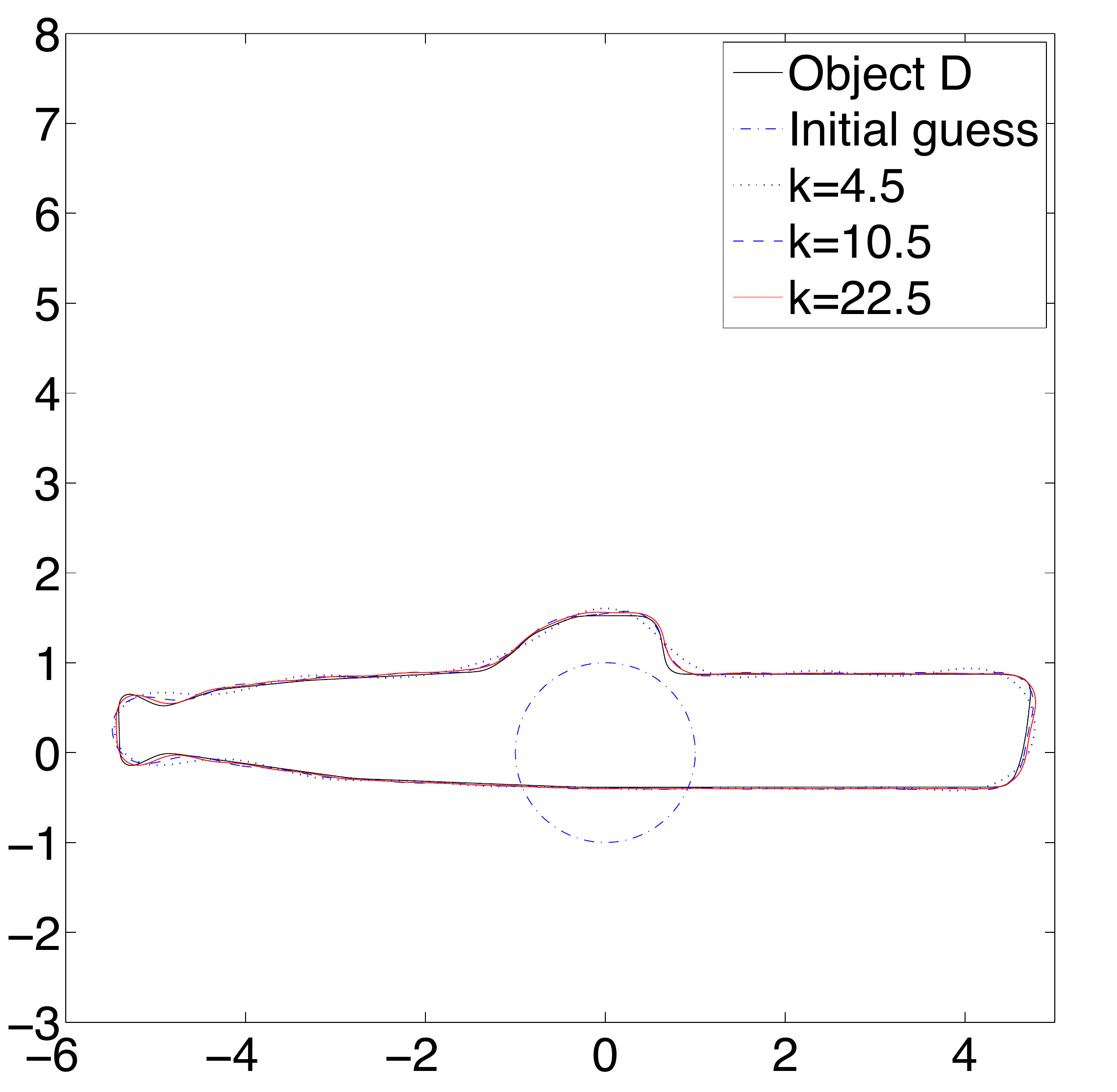}\label{sub1a}
      }
  \subfigure[Back part of the submarine]{      
            \includegraphics[scale=0.3]{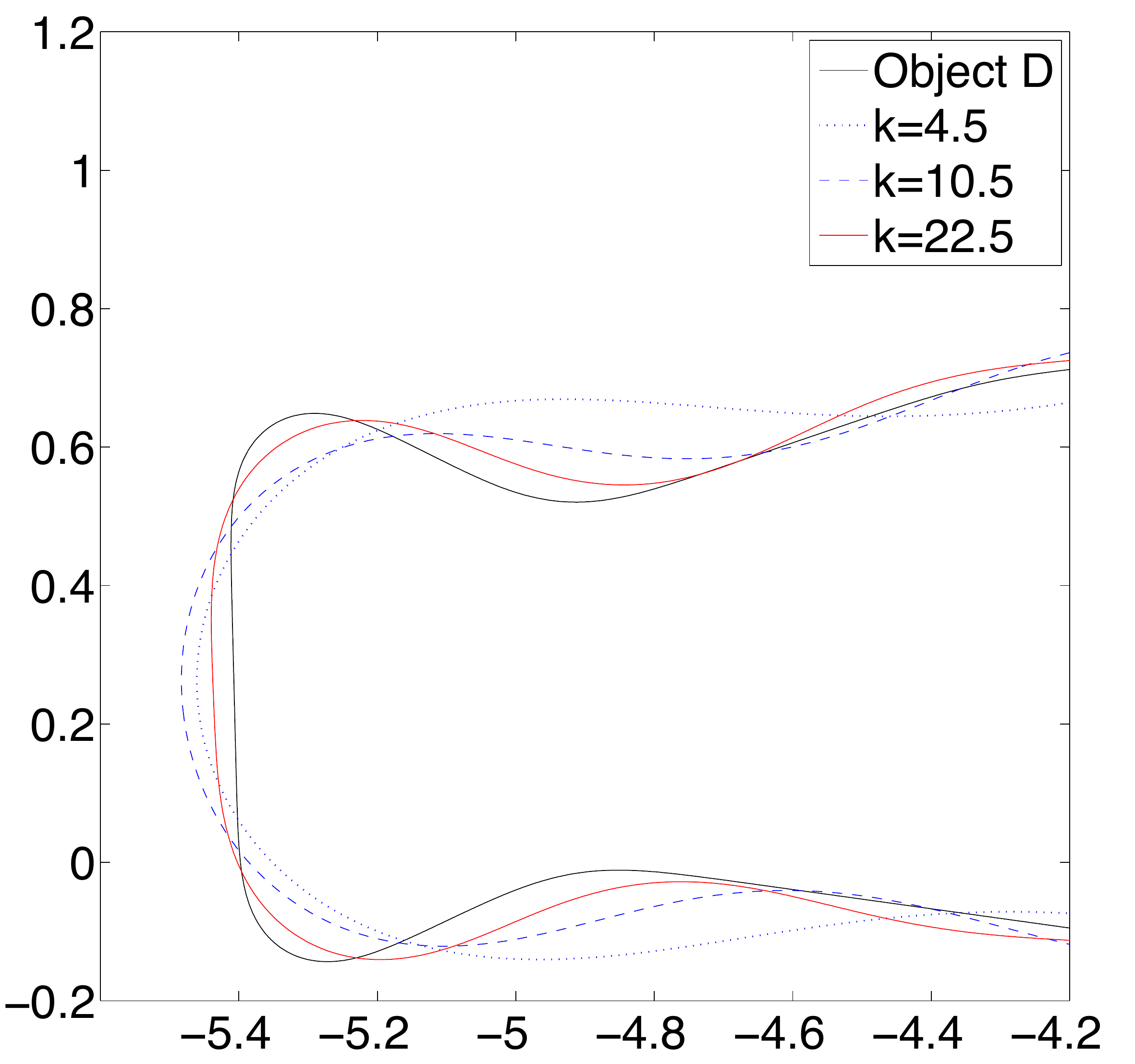}\label{sub1b}
  }
    \subfigure[Front part of the intersection of the tower and hull of the submarine]{      
            \includegraphics[scale=0.3]{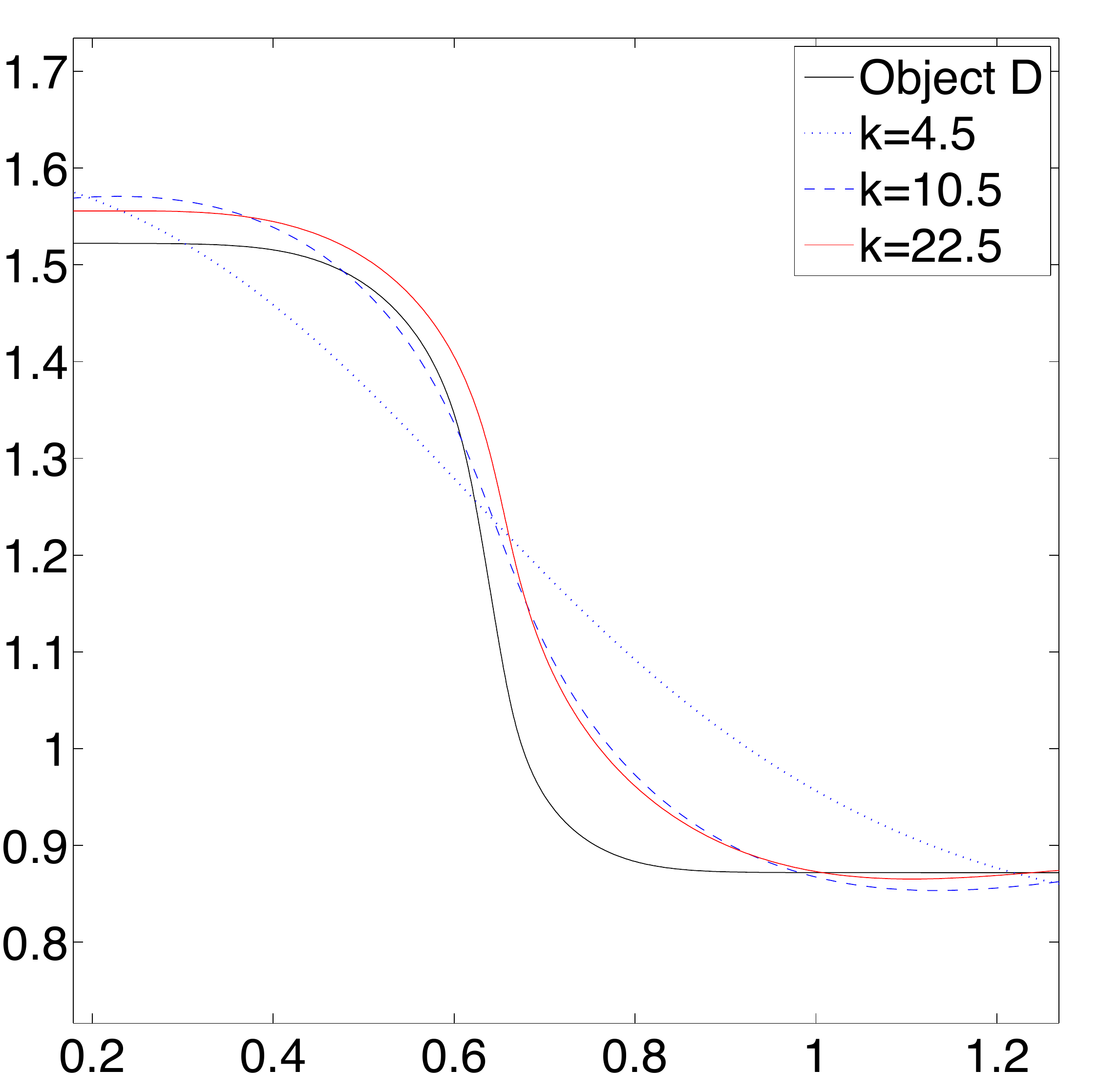}\label{sub1c}
  }
      \subfigure[Back part of the intersection of the tower and hull of the submarine]{      
            \includegraphics[scale=0.3]{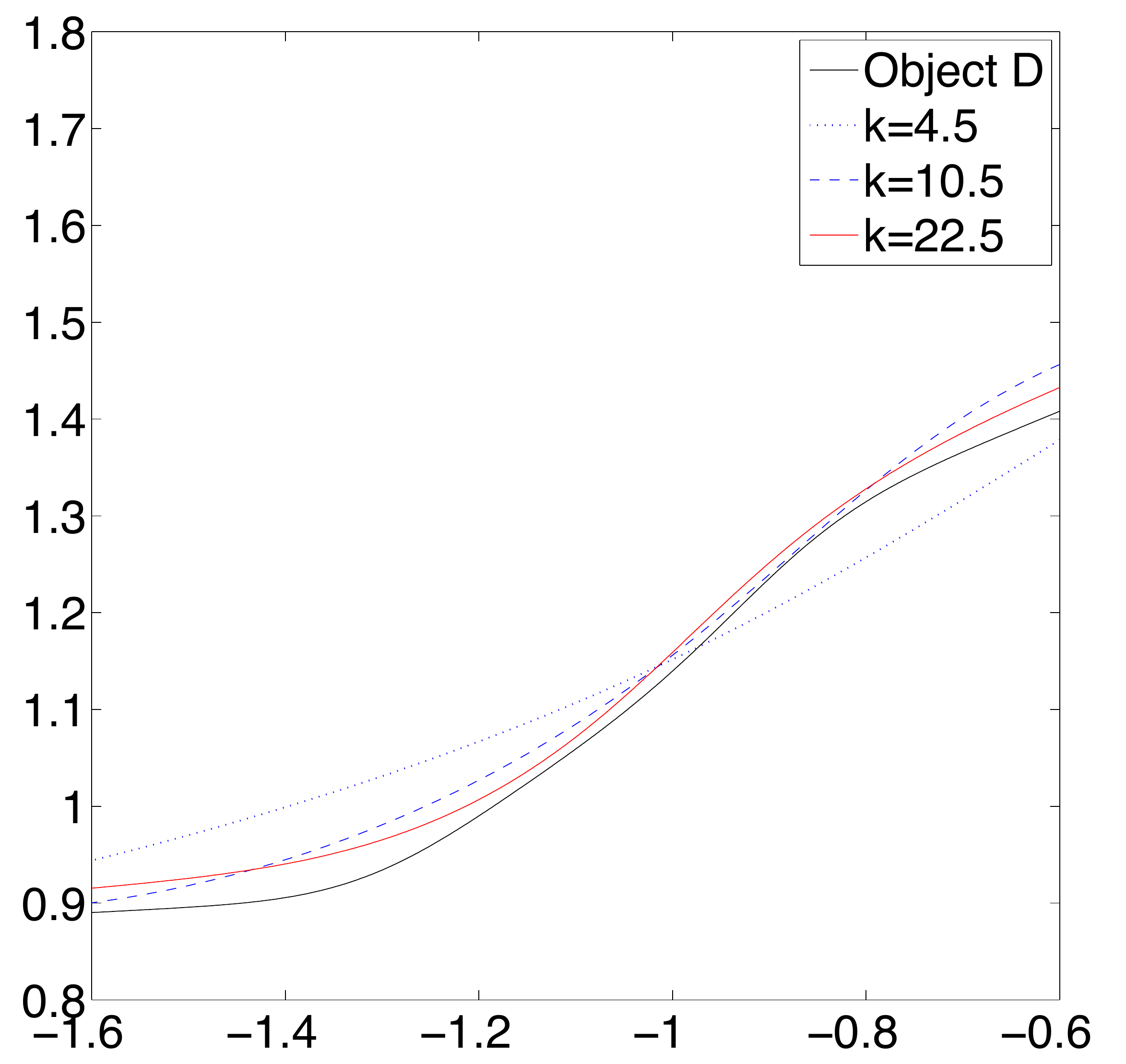}\label{sub1d}
  }
  \caption[Applying the method to reconstruct the submarine.]{\footnotesize 
Reconstruction using recursive linearization with Newton's method 
to reconstruct the submarine-like object.} \label{sub}
\end{figure}

In our final example,
we reconstruct an object with a submarine-like shape. 
Our initial guess is the circle of 
radius 1 centered at the origin. 
The initial wavenumber is $k_0=0.1$, with $\Delta k=0.4$, and $J=57$. 
We let the number of far field measurements be $M=128$ and assumed
six directions of incidence ($L=6$). 
We choose the damping parameter $\rho=1$.  
We use the bandlimit $b=2\ceil{k}+1$ at each wavenumber, with a smooth roll-off in
in the curve resammpling filter from $b$ to $b+50$. 
Within Newton's method, the stopping criterion used 
for the nonlinear residual is
$\epsilon = \| F(\Gamma) - u^\infty \| < 10^{-3}$.
We used $N_1=\ceil{20k|\Gamma_j|}$ quadrature points in solving
the integral equations within the Newton step. 
In Fig. \ref{sub1a}, 
we plot the solution obtained for various maiximal frequencies. 
Figs.  \ref{sub1b}--\ref{sub1d} show details of the reconstruction 
in certain high curvature regions of the object.

\begin{figure}[h!]
  \centering
    \subfigure[Number of iterations for the aircraft-like object]{      
            \includegraphics[scale=0.28]{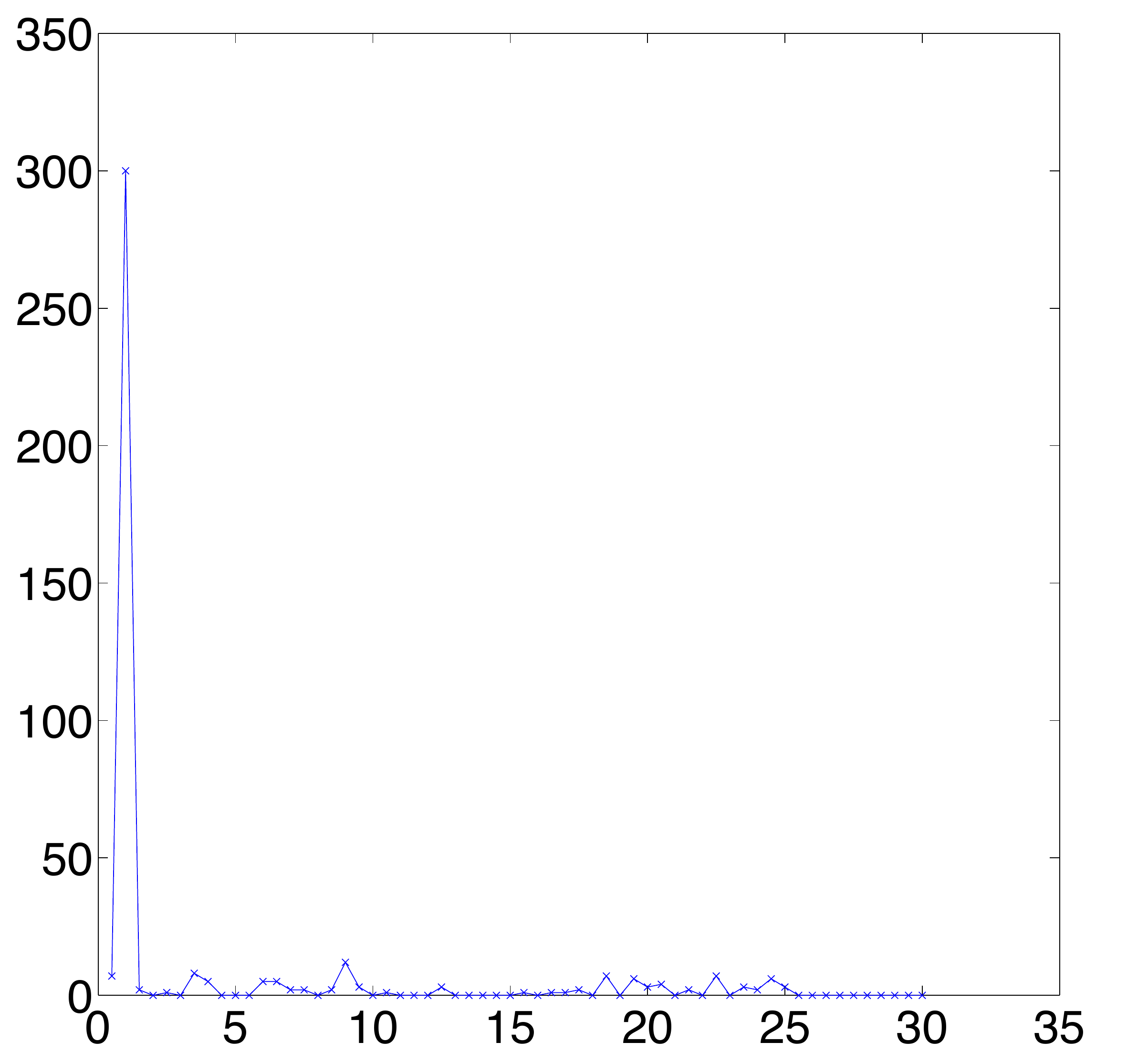}\label{planeit}
  }
\subfigure[Number of iterations for the submarine-like object]{
      \includegraphics[scale=0.28]{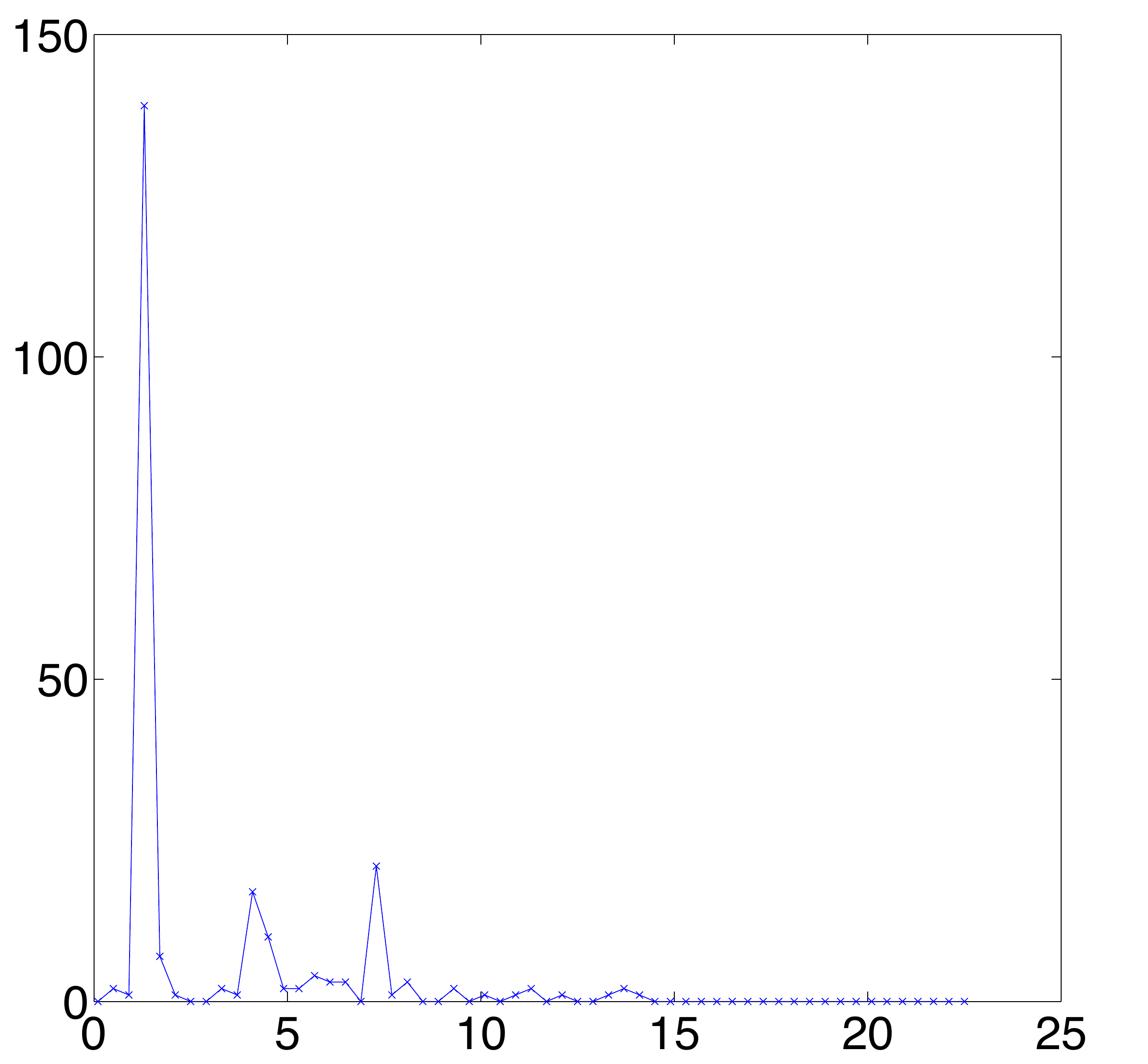}\label{subit}
      }
  \caption[Number of iterations for the Newton's method at each wavenumber $k$.]{\footnotesize 
Number of iterations for the Newton's method at each wavenumber $k$. } \label{numberit}
\end{figure}

The number of Newton iterations at each frequency to recover the aircraft-like and 
submarine-like objects are shown in Fig. \ref{numberit}.

Note that a large number of damped steps are required at the lowest frequency, when the initial guess is very far 
from the correct shape, but that the number decreases significantly after that.

\section{Conclusion}
We have presented a technique for reconstructing the shape of 
two-dimensional sound-soft obstacles, given the far-field pattern, using
multiple angles of incidence and multiple frequencies.
While the problem is both 
ill-posed and nonlinear, a combination of techniques makes it tractable.
First, we reduce the number of degrees of freedom to be determined based
on physical considerations (the approximate dimensions of the object
in wavelengths). We do not make use of Tikhonov regularization or 
other generic regularization schemes. Second, when multiple frequency data is available, 
we make use of recursive linearization \cite{Chen}. 
In our experiments, Newton's method with damping has trouble converging at
the lowest frequency, when the initial
guess is far from the desired minimum. Subsequently, recursive 
linearization enables rapid convergence, consistent with the analysis
of Sini and Th\`{a}nh \cite{Sini11}.
In our method, we also make use of the fact that the conditioning of the reconstruction
problem is improved by using multiple angles of incidence. 
Finally, we employ high-order accurate discretizations and the HODLR fast solver
to accelerate the solution of the necessary forward scattering problems, 
so that all our experiments are easily carried out using modest computational resources.
Unlike many algorithms for the inverse obstacle scattering problem, we make no
assumption about the parametrization of the geometry, such as being star-shaped. 
In our present formulation,
however, we do assume that the obstacle is a single, simply-connected region in the plane.

The method is easily extended to the case of sound-hard obstacles or 
impedance boundary conditions. Moreover, most aspects have straightforward
three-dimensional analogs: recursive linearization, high-order discretization and 
fast, forward scattering solvers (although the latter two are still areas of active 
research). One aspect that is not so straightforward concerns
the development of a simple and robust parametrization of surfaces
which permits ``physical regularization" - systematically reducing the 
number of degrees of freedom to be solved for in each linearized problem by
enforcing some kind of band-limit on the perturbation.
We are currently considering a variety of 
approaches for this and progress will be reported at a later date.

Two issues we have not addressed here are inverse obstacle scattering when only partial
aperture data is available and when only the magnitude of the far field is measured,
rather than magnitude and phase. These are of significant practical importance and also
active areas of research.
 \label{sec:concl}


%
%

\section*{Acknowledgments}
The work reported was supported by the Air Force Office of Scientific Research under grant number FXXXXX-XX- X-XXXX. The authors also gratefully acknowledge the valuable comments and suggestions of Michael O'Neil and Sivaram Ambikasaran.
\bibliography{./Bib}

\end{document}